\documentclass[11pt]{amsart}
\usepackage{amsmath} \usepackage{amsfonts,amssymb}
\usepackage{enumerate} \usepackage{amsthm}
\theoremstyle{plain} \newtheorem{theorem}{Theorem}[section]
\newtheorem{proposition}[theorem]{Proposition}
 \newtheorem{lemma}[theorem]{Lemma}
\newtheorem{corollary}[theorem]{Corollary}
\newtheorem{definition}[theorem]{Definition}

\theoremstyle{definition} \newtheorem{assu}[theorem]{Assumption}
\newtheorem{nota}[theorem]{Notation}
\newtheorem{remark}[theorem]{Remark}
\numberwithin{equation}{section}

\DeclareMathOperator{\diff}{d}
\DeclareMathOperator{\cn}{\mathrm{div}}
\DeclareMathOperator{\curl}{\mathrm{curl}}
\DeclareMathOperator{\Op}{\mathrm{Op}} \DeclareMathOperator{\HF}{HF}
 \DeclareMathOperator{\LF}{LF}

\newcommand{\Cr}{\mathcal{C}}

 \newcommand{\Kr}{\mathcal{K}}

 \newcommand{\Rr}{\mathcal{R}}
\newcommand{\Sr}{\mathcal{S}} 
 
 \newcommand{\Xr}{\mathcal{X}}
 \newcommand{\Zr}{\mathcal{Z}}
\def\bi{{\chi}}

 \def\cpe{\nu} 
\def\defe{=}
\def\defn{\mathrel{:=}}

\def\Di{F} \def\DII{G}

\def\eps{\varepsilon} \def\ffp#1{\partial_{t}{#1}+v\cdot\nabla {#1}}

\def\fp#1#2{\frac{\partial{#1}}{\partial #2}}
\def\fpll#1#2{(\partial{#1}/\partial #2)}
\def\fplll#1#2{\partial{#1}/\partial #2}
\def\Fi#1#2{\Lambda^{{#2}}_{{#1}}}

\def\gi{{g}} \def\ge{\geqslant}
 \def\id{I}
\def\ik{m} 
 \def\ki{{k}}
\def\la{\left\lvert}
\def\lA{\left\lVert} 
\def\le{\leqslant}
\def\les{\lesssim}
\def\L#1{\langle{#1}\rangle}

  \def\nh#1{{H^{{#1}}_{~}}} 
\def\nhsc#1{{{H}^{{#1}}_{\cpe}}}
\def\norm#1{\left\lVert{\smash[t]{{#1}}}\right\rVert}

\def\pa{a}  \def\pe{\kappa} \def\PA{A} \def\r{\mu}
\def\ra{\right\rvert} \def\rA{\right\rVert} \def\scal#1#2{\langle \,
{#1} \hspace{1pt},\hspace{1pt} {#2} \,\rangle} \def\te{{\eta}}
 \def\tvar{\widetilde{\var}}
\def\tvari{\widetilde{p}} \def\tvard{\widetilde{v}} \def\tvariii{\widetilde{\theta}}
   \def\tvar{\widetilde{\var}}
 \def\tz{{\zeta}} \def\var{U} \def\vari{p}
\def\variii{\theta} \def\vard{v} \def\virgp{\raise 2pt\hbox{,}}
  \def\xN{\mathbb{N}}
 \def\xR{\mathbb{R}} 
\def\xT{\mathbb{T}}

\begin{document}
\title[Low Mach number flows]{Low Mach number flows, and combustion}
\author[T. Alazard]{Thomas Alazard}
\address{MAB, Universit\'e de Bordeaux I\\ $33405$ Talence Cedex,
France} \email{thomas.alazard@math.u-bordeaux1.fr}

\begin{abstract}
We prove uniform existence results for the full Navier-Stokes
equations for 
time intervals which are independent 
of the Mach number, the Reynolds number and the P\'eclet number. 
We consider general equations of state and 
we give an application for 
the low Mach number limit combustion problem introduced by Majda in~\cite{Maj}.
\end{abstract}

\maketitle

\section{Introduction}
%
%
For a fluid with density $\varrho$, velocity $v$, pressure $P$,
temperature $T$, internal energy $e$, Lam\'e coefficients $\zeta,\eta$
and coefficient of thermal conductivity~$k$, the full Navier-Stokes
equations, written in a non-dimensional way, are
\begin{equation}\label{system:ANS}
\left\{
\begin{aligned} 
&\partial_{t}\rho+\cn(\rho v)=0,\\ &\partial_{t}(\rho v)+\cn(\rho
v\otimes v)+ \frac{\nabla P}{\eps^{2}}=\mu\bigl(2\cn(\zeta
Dv)+\nabla(\eta\cn v)\bigr),\\ &\partial_{t}(\rho e)+\cn(\rho v
e)+P\cn v=\kappa\cn(k\nabla T)+ Q,
\end{aligned}
\right.
\end{equation}
where $\eps\in (0,1]$, $(\mu,\kappa)\in [0,1]^{2}$ and $ Q$ is a given
source term (see \cite{Embid1,Klein,Maj}). In order to be closed, the
system is supplemented with two equations of state, so that
$\rho,P,e,T$ are completely determined by only two of these
variables. Also, it is assumed that $\zeta$, $\eta$ and $k$ are smooth
functions of the temperature.

\smallbreak 

This paper is devoted to the asymptotic limit where the
Mach number $\eps$ tends to $0$. We are interested in proving results
independent of the Reynolds number $1/\mu$ and the P\'eclet
number~$1/\kappa$. Our main result asserts that the classical
solutions of \eqref{system:ANS} exist and are uniformly bounded on a
time interval independent of $\eps$, $\mu$ and~$\kappa$.

This is a continuation of our previous work~\cite{TA2} where the study
was restricted to perfect gases and small source terms $Q$ of size
$O(\eps)$.  We refer to the introduction of~\cite{TA2} for references
and a short historical survey of the background of these problems (see
also the survey papers of Danchin \cite{Danchin}, Desjardins and Lin~\cite{DL}, 
Gallagher~\cite{Gallagher}, Schochet~\cite{Schochet} and Villani~\cite{Villani}).

\smallbreak

The case of perfect gases is interesting in its own: first, perfect
gases are widely studied in the physical literature; and second, it
contains the important analysis of the singular terms. Yet, modeling
real gases requires general equations of state (see~\cite{BLP,MY}). Moreover, we shall 
see that it is interesting to consider large source
terms~$Q$ for it allows us to answer a question addressed by Majda in~\cite{Maj} 
concerning the combustion equations. 

\subsection{The equations}

To be more precise, we begin by rewriting the equations under the form
$
L(u,\partial_{t},\partial_{x})u+\eps^{-1}S(u,\partial_{x})u=0,
$
which is the classical framework of a singular limit problem. 

%
%
Before we proceed, three observations are in order. Firstly, 
for the low Mach number limit 
problem, the point is not so much to use the conservative form of the
equations, but instead to balance the acoustics components. This is
one reason it is interesting to work with the unknowns $P,v$,~$T$ (see~\cite{Maj}). 
Secondly, the general case must allow for large density and temperature
variations as well as very large acceleration of
order of the inverse of the Mach number (see Section~$5$ in~\cite{Klein}). Since
$\partial_{t}v$ is of order of $\eps^{-2}\nabla P$, this suggests that
we seek $P$ under the form $P={\rm Cte}+O(\eps)$. As in~\cite{MS1}, 
since $P$ and $T$ are positive functions, it is
pleasant to set 
\begin{equation}\label{defi:ptheta}
P=\underline{P}e^{\eps p},\quad T=\underline{T}e^{\theta},
\end{equation}
where $\underline{P}$ and $\underline{T}$ are given positive
constants, say the reference states at spatial infinity. 
Finally, the details of the following computations are given in the Appendix.
%
%

From now on, the unknown is $(p,v,\theta)$ with values in~$\xR\times\xR^{d}\times\xR$.
We are interested in the general case where $p$ and $\theta$ are uniformly
bounded in $\eps$ (so that $\nabla T=O(1)$ and
$\partial_{t}v=O(\eps^{-1})$). 
 
\smallbreak

By assuming that $\rho$ and $e$ are given smooth functions of $(P,T)$, 
it is found that, for smooth solutions of~\eqref{system:ANS}, 
$(P,v,T)$ satisfies a system of the form:
\begin{equation}\label{system:ANSF}
\left\{
\begin{aligned} 
&\alpha(\partial_{t}P+v\cdot\nabla P)+\cn v
=\kappa\beta\cn(k\nabla T)+\beta Q,\\
&\rho(\partial_{t}v+v\cdot\nabla v)+\frac{\nabla
  P}{\eps^{2}}=\mu\bigl(2\cn(\zeta Dv)+\nabla(\eta\cn v)\bigr),\\
&\gamma(\partial_{t}{T}+v\cdot\nabla T)+\cn v=\kappa\delta\cn(k\nabla T)+\delta Q,
\end{aligned}
\right.
\end{equation}
where the coefficients $\alpha$, $\beta$, $\gamma$ and $\delta$ are smooth functions of $(P,T)$. 
Then, by writing $\partial_{t,x}P=\eps
P\partial_{t,x} p$, $\partial_{t,x}T=T\partial_{t,x}\theta$ 
and redefining the functions $k$, $\zeta$ and $\eta$, 
it is found that $(p,v,\theta)$ satisfies a system of the form:
\begin{equation}\label{system:NSint}
\left\{
\begin{aligned}
&g_{1}(\phi)(\partial_{t}p+v\cdot\nabla p)
+\frac{1}{\eps}\cn v =
\frac{\kappa}{\eps}\chi_{1}(\phi)\cn(k(\theta)\nabla\theta)+\frac{1}{\eps}\chi_{1}(\phi) Q,\\
&g_{2}(\phi)(\partial_{t}v +v\cdot\nabla v) +\frac{1}{\eps}\nabla p 
= \mu B_{2}(\phi,\partial_{x})v,\\
&g_{3}(\phi)(\partial_{t}\theta +v\cdot\nabla \theta)+\cn v
=\kappa \chi_{3}(\phi)\cn(k(\theta)\nabla\theta)+\chi_{3}(\phi) Q,
\end{aligned}
\right.
\end{equation}
where $\phi\defn (\theta,\eps p)$ and 
$
B_{2}(\phi,\partial_{x})=\chi_{2}(\phi)\cn (\zeta(\theta) D \cdot)
+\chi_{2}(\phi)\nabla(\eta(\theta) \cn \cdot).
$

\smallbreak

We are now in position to explain the main differences between ideal gases and 
general gases. 
Firstly, we note that the source term $Q$ 
introduces an arbitrary unsigned large term 
of order of $1/\eps$ in the equations. 
Secondly, to emphasize the role of the thermodynamics, we suppose now
that $Q=0$ and we mention that, for perfect gases, the coefficient
$\chi_{1}(\phi)$ is a function of $\eps p$ alone (see Proposition~\ref{prop:structurepg}). 
Hence, for perfect gases, 
the limit constraint is linear in the sense that it 
reads $\cn v_{e}=0$ with $v_{e}\defe v-\kappa
\chi_{1}(0)k(\theta)\nabla\theta$. By contrast, for general equations of
state, the limit constraint is nonlinear.

\subsection{Assumptions}\label{assu:assus}

To avoid confusion, we denote by $(\vartheta,\wp)\in\xR^{2}$ the 
place holder of the unknown $(\theta,\eps p)$. Hereafter, it is assumed that:

\begin{enumerate}[$({{\rm H}}1)$]
\item 
The functions $\tz$, $\te$ and $\ki$ 
are $C^{\infty}$ functions of $\vartheta\in\xR$, satisfying 
$\ki>0$, $\tz>0$ and $\te+2\tz>0$. 

\item \label{assus}
The functions $\gi_{i}$ and $\bi_{i}$ ($i=1,2,3$) are $C^{\infty}$ 
positive functions of $(\vartheta,\wp)\in\xR^{2}$. Moreover, 
$$
\bi_{1}<\bi_{3},
$$
and there exist two functions $\Di$ and $\DII$ such that  
$(\vartheta,\wp)\mapsto(\Di(\vartheta,\wp),\wp)$ and 
$(\vartheta,\wp)\mapsto(\vartheta,\DII(\vartheta,\wp))$ 
are $C^{\infty}$ diffeomorphisms from $\xR^{2}$ onto $\xR^{2}$, $F(0,0)=G(0,0)=0$ and
\begin{equation*}
\gi_{1}\frac{\partial \Di}{\partial \vartheta}
=-\gi_{3}\frac{\partial \Di}{\partial \wp} >0,\qquad 
\gi_{1}\bi_{3}\frac{\partial \DII}{\partial \vartheta}
=-\gi_{3}\bi_{1}\frac{\partial \DII}{\partial \wp} <0.
\end{equation*}
\end{enumerate}

\begin{remark}
Assumption~(H\ref{assus}) is used to prove various energy estimates. 
The main hypothesis is the inequality
$\bi_{1}<\bi_{3}$. In Appendix~\ref{appendix:COV}, it is proved that the inequality $\bi_{1}<\bi_{3}$ 
holds whenever 
the density $\rho$ and the energy $e$ are $C^{\infty}$ functions 
of $(P,{T})\in (0,+\infty)^{2}$, such that 
$\rho>0$ and 
\begin{equation}\label{P1FI}
P\fp{\rho}{P}+T\fp{\rho}{T}=\rho^{2}\fp{e}{P}\virgp
\quad\fp{\rho}{P} > 0, \quad \fp{\rho}{{T}} < 0,
\quad
\fp{e}{{T}}\fp{\rho}{P}
>\fp{e}{P}\fp{\rho}{{T}}\cdot
\end{equation}
\end{remark}

\subsection{Main result}
We are interested in the case without smallness assumption: 
namely, we consider general initial data, general equations of state and large source terms $Q$. 
To get around the above mentioned  nonlinear features of the penalization operator, 
we establish a few new qualitative properties. These properties are enclosed 
in various uniform stability results, which assert that the
classical solutions of~\eqref{system:NSint} exist 
and they are uniformly bounded for a time independent
of~$\eps$,~$\mu$ and~$\kappa$. We concentrate 
below on the whole space problem or the periodic case and 
we work in the Sobolev spaces $H^{\sigma}$ 
equipped with the norms $\lA u\rA_{H^{\sigma}}
\defn \lA (\id-\Delta)^{\sigma/2}u\rA_{L^{2}}$.

The following result is the core of all our other uniform stability results. 
On the technical side, it contains the idea 
that one can prove uniform estimates without 
uniform control of the $L^{2}_{x}$ norm of the velocity~$v$. 

\begin{theorem}\label{theo:main}
Let $d=1$ or $d\ge 3$ and $\xN\ni s>1+d/2$. 
For all source term $ Q= Q(t,x)\in C^{\infty}_{0}(\xR\times\xR^{d})$ 
and all $M_{0}>0$, there exist $T>0$ and $M>0$ 
such that, for all $(\eps,\mu,\kappa)\in (0,1]\times [0,1]\times [0,1]$ and all initial data 
$(p_{0},v_{0},\theta_{0})\in H^{s+1}(\xR^{d})$ satisfying
\begin{equation}\label{TM:CI}
\norm{(\nabla p_{0},\nabla v_{0})}_{H^{s-1}} 
+ \norm{(\theta_{0},\eps p_{0},\eps v_{0})}_{H^{s+1}}\le M_{0},
\end{equation}
the Cauchy problem for~\eqref{system:NSint} has a unique classical solution  
$(p,v,\theta)$ in
$C^{0}([0,T];H^{s+1}(\xR^{d}))$  such that
\begin{equation}\label{TH:NORM}
\sup_{t\in [0,T]} \norm{(\nabla p(t),\nabla v(t))}_{H^{s-1}} 
+ \norm{(\theta(t),\eps p(t),\eps v(t))}_{H^{s}}\le M.
\end{equation}
\end{theorem}
A refined statement is proved in Section~\ref{section:proof}. 

A notable corollary of Theorem~\ref{theo:main} is Theorem~\ref{theo:GP}, 
which is the requested result for application to the 
low Mach number limit. 
Detailed discussions of the periodic case and the combustion equations
are included in Sections~\ref{section:Torus} and~\ref{section:Combustion}. 
The assumption $d\neq 2$ is explained in Remark~\ref{R2D}.

\section{Preliminaries}

In order not to interrupt the proofs later on, 
we collect here some estimates. 
The main result of this section is Proposition~\ref{prop:A2F}, which 
complements the Friedrichs-type estimate
\begin{equation}\label{UEF}
\lA \nabla v\rA_{H^{s}}\le \lA \cn v\rA_{H^{s}}+\lA \curl v\rA_{H^{s}},
\end{equation}
which is immediate using Fourier transform. 
We prove a variant where $\cn v$ 
is replaced by $\cn(\rho v)$ 
where $\rho$ is a positive weight. 

\smallbreak
\noindent \textbf{Notation.} 
The symbol $\les$ stands for $\le$ up to a positive, multiplicative constant, 
which depends only on parameters that are considered fixed.

\subsection{Nonlinear estimates}
Throughout the paper, 
we will make intensive and often implicit uses of the following estimates. 

For all $\sigma\ge 0$, there exists $K$ such that, 
for all $u,v\in L^{\infty}\cap H^{\sigma}(\xR^{d})$,
\begin{equation}\label{PTE}
\lA u v\rA_{H^{\sigma}}\le K\lA u\rA_{L^{\infty}}\lA v\rA_{H^{\sigma}}
+K\lA u\rA_{H^{\sigma}}\lA v\rA_{L^{\infty}}.
\end{equation}

For all $s>d/2$, $\sigma_{1}\ge 0$, $\sigma_{2}\ge 0$ such that $\sigma_{1}+\sigma_{2}\le 2s$, 
there exists a constant $K$ such that, 
for all $u\in H^{s-\sigma_{1}}(\xR^{d})$ and $v\in H^{s-\sigma_{2}}(\xR^{d})$,
\begin{equation}\label{PRS}
\lA u v\rA_{H^{s-\sigma_{1}-\sigma_{2}}}\le K 
\lA u\rA_{H^{s-\sigma_{1}}}
\lA v\rA_{H^{s-\sigma_{2}}}.
\end{equation}

For all $s>d/2$ and for all $C^{\infty}$ function $F$ vanishing at the origin, 
there exists a smooth function $C_{F}$ such that, for all $u\in H^{s}(\xR^{d})$,
\begin{equation}\label{PTE2}
\lA F(u)\rA_{H^{s}}\le C_{F}(\lA u\rA_{L^{\infty}})\lA u\rA_{H^{s}}.
\end{equation}

\subsection{Estimates in~$\xR^{3}$}
Consider the Fourier multiplier $\nabla\Delta^{-1}$ with symbol $-i\xi/\la
\xi\ra^{2}$. This operator is, at least formally, 
a right inverse for the divergence operator. The only think we will use below 
is that $\nabla\Delta^{-1}u$ 
is well defined whenever 
$u=u_{1}u_{2}$ with $u_{1},u_{2}\in L^{\infty}\cap H^{\sigma}(\xR^{d})$ for some 
$\sigma\ge 0$.

\begin{proposition}\label{prop:RD}
Given $d\ge 3$ and $\sigma\in\xR$, the Fourier 
multiplier $\nabla\Delta^{-1}$ is well defined on $L^{1}(\xR^{d})\cap
H^{\sigma}(\xR^{d})$ with values in $H^{\sigma+1}(\xR^{d})$. 
Moreover, there exists a constant $K$ such that, for
all $u\in L^{1}(\xR^{d})\cap H^{\sigma}(\xR^{d})$,
\begin{equation}\label{propRD1}
\lA \nabla\Delta^{-1} u\rA_{H^{\sigma+1}}\le K\lA
u\rA_{L^{1}}+K\lA u\rA_{H^{\sigma}}.
\end{equation}
\end{proposition}
\begin{proof}[Proof]
Set $\L{\xi}\defn (1+\la \xi\ra^{2})^{1/2}$. It suffices 
to check that the $L^{2}$-norm of 
$(\L{\xi}^{\sigma+1}/\la \xi\ra)\la \widehat{u}(\xi)\ra$ 
is estimated by the right-hand side of \eqref{propRD1}. 
To do that we write
$$
\int_{\la \xi\ra\le 1}\frac{\L{\xi}^{2\sigma+2}}{\la \xi\ra^{2}}
\la \widehat{u}(\xi)\ra^{2}\,d\xi\les\lA u\rA_{L^{1}}^{2},
\quad
\int_{\la\xi\ra\ge 1}\frac{\L{\xi}^{2\sigma+2}}{\la \xi\ra^{2}}
\la \widehat{u}(\xi)\ra^{2}\,d\xi\les\lA u\rA_{H^{\sigma}}^{2},
$$
where we used $1/\la \xi\ra^{2} \in L^{1}(\{\la\xi\ra\le 1\})$  for all $d\ge 3$.
\end{proof}

The next proposition is well known. Its corollary is a 
special case of a general estimate established in \cite{BDD2}. 
\begin{proposition}\label{prop:Linftycontrol}
Given~$d\ge 3$ and~$s>d/2$, there exists a constant 
$K$ such that, for all~$u\in H^{s}(\xR^{d})$, 
\begin{equation}\label{restriction:dimge3}
\norm{u}_{L^{\infty}} \le K \norm{\nabla u}_{\nh{s-1}}.
\end{equation}
\end{proposition}
\begin{proof}[Proof]Since $H^{s}(\xR^{d})\hookrightarrow 
L^{\infty}(\xR^{d})$, it suffices to prove the result for~$u$ 
in the Schwartz class~$\Sr(\xR^{d})$. Now, starting from the Fourier inversion theorem, 
the Cauchy--Schwarz inequality yields the desired estimate:
\begin{align*}
\lA u \rA_{L^{\infty}}
\le \left( \int \frac{d\xi}{\la\xi\ra^{2}\L{\xi}^{2(s-1)}}\right)^{1/2}
\left(\int\L{\xi}^{2(s-1)}\la\xi\widehat{u}(\xi)\ra^{2}\,d\xi\right)^{1/2}
\les \lA \nabla u \rA_{H^{s-1}}.
\end{align*}
\renewcommand{\qed}{}
\end{proof}
\begin{corollary}
Given $d\ge 3$ and $\xN\ni s>d/2$, there exists a constant $K$ such that,
for all $u_1,u_2\in H^{s}(\xR^{d})$,
\begin{equation}\label{inter}
\lA u_{1}u_{2}\rA_{H^{s}}\le K\lA \nabla u_{1}\rA_{\nh{s-1}}
\lA u_{2} \rA_{H^{s}}.
\end{equation}
\end{corollary}
\begin{proof}[Proof]
One has to estimate the $L^{2}$-norm of $\partial_{x}^{\alpha}(u_1 u_2)$, 
where $\alpha\in\xN^{d}$ satisfies $\la \alpha\ra\le s$. Rewrite this term as
$u_1 \partial_{x}^{\alpha} u_2 + [\partial_{x}^{\alpha},u_1 ] u_2$. 
Since the commutator is a sum of terms of the form $\partial_{x}^{\beta}u_{1}\partial_{x}^{\gamma}u_{2}$ 
with $\beta>0$, the product rule~\eqref{PRS} implies that
\begin{equation}\label{Kato_Ponce}
\lA  [\partial_{x}^{\alpha},u_1] u_2 \rA_{L^{2}}\les
\lA \nabla u_1\rA_{H^{s-1}}\lA u_2\rA_{H^{s}}.
\end{equation}
Moving to the estimate of the first term, we write
\begin{equation*}
\lVert u_{1}\partial_{x}^{\alpha}u_{2}\rVert_{L^{2}}\le \lA
u_{1}\rA_{L^{\infty}}\lVert u_{2}\rVert_{H^{s}}
\les \lA \nabla u_1 \rA_{H^{s-1}}\lVert
u_{2}\rVert_{H^{s}}.\qquad\qquad \square
\end{equation*}
\renewcommand{\qed}{}
\end{proof}

\subsection{A Friedrichs' Lemma}
With these preliminaries established, we are prepared to prove the following:
\begin{proposition}\label{prop:A2F}
Let $d\ge 3$ and $\xN\ni s>d/2$. There exists a function $\Cr$ such that, 
for all 
$\varphi\in H^{s+1}(\xR^{d})$ and all vector field $v\in H^{s+1}(\xR^{d})$, 
\begin{align}
\lA \nabla v\rA_{H^{s}}\le C \lA \cn(e^{\varphi} v)\rA_{H^{s}}
+C\lA \curl v\rA_{H^{s}},\label{A2F1}
\end{align} 
where 
$C\defn (1+\lA \varphi\rA_{H^{s+1}})
\Cr\bigl(\lA \varphi\rA_{H^{s}},\lA\nabla\varphi\rA_{L^{\infty}}\bigr)$.
\end{proposition}
\begin{proof}
For this proof, we use the notation
$$
R\defe \lA \cn(e^{\varphi} v)\rA_{H^{s}}
+\lA \curl v\rA_{H^{s}},
$$
and we denote by $C_{\varphi}$ various constants depending only on
$\lA \varphi\rA_{H^{s}}+\lA\nabla\varphi\rA_{L^{\infty}}$.

All the computations 
given below are meaningful since it is sufficient to prove \eqref{A2F1} 
for $C^{\infty}$ functions with compact supports. We begin by setting
$$
\widetilde{v}\defe v +\nabla\Delta^{-1}\bigl(\nabla\varphi\cdot v\bigr).
$$
The reason to introduce $\widetilde{v}$ is that
\begin{equation*}
e^{\varphi} \cn \widetilde{v} =\cn (e^{\varphi} v), \qquad  
\curl \widetilde{v}=\curl v.\label{A2F1E1}
\end{equation*}
Hence, by using~\eqref{UEF}, we have
\begin{equation}\label{F.1}
\lA \nabla \widetilde{v}\rA_{H^{s}}\le 
\lA e^{-\varphi}\cn (e^{\varphi}
v)\rA_{H^{s}}+\lA \curl v\rA_{H^{s}}\le C_{\varphi}R.
\end{equation}
The proof of \eqref{A2F1} thus reduces to estimating $v_{1}\defn v-\widetilde{v}$, which satisfies
\begin{equation*}
\cn (e^{\varphi} v_{1})=-e^{\varphi}\nabla \varphi\cdot\widetilde{v}, 
\qquad \curl v_{1}=0.\label{A2F1E2}
\end{equation*}
Again, to estimate $v_{1}$ we introduce $\widetilde{v}_{1}\defn
v_{1}+\nabla\Delta^{-1}\bigl(\nabla\varphi\cdot v_{1}\bigr)$, which solves
$$
\cn \widetilde{v}_{1} =-\nabla\varphi\cdot\widetilde{v},  
\qquad \curl \widetilde{v}_{1} =0.
$$
The estimate \eqref{UEF} implies that 
$\lA\nabla \widetilde{v}_{1}\rA_{H^{s}}\le 
\lA \nabla \varphi\cdot\widetilde{v}\rA_{H^{s}}$. 
By using \eqref{F.1} and the 
product rule \eqref{inter}, applied with $u_{1}=\widetilde{v}$ and $u_{2}=\nabla\varphi$, 
we find that 
\begin{equation}\label{A2F1PR1}
\lA \nabla\widetilde{v}_{1}\rA_{H^{s}}
\le \lA \nabla \varphi\cdot\widetilde{v}\rA_{H^{s}}
\les\lA\varphi\rA_{H^{s+1}}
\lA\nabla\widetilde{v}\rA_{H^{s-1}}\le \lA\varphi\rA_{H^{s+1}}C_{\varphi}R. 
\end{equation}
Hence, it remains only to estimate 
$v_{2}\defe v_{1}-\widetilde{v}_{1}$, which satisfies
$$
\cn (e^{\varphi} v_{2})=-e^{\varphi}\nabla\varphi\cdot \widetilde{v_{1}} \quad\text{and} \quad \curl v_{2}=0.
$$
To estimate $v_{2}$ the key point is the estimate
\begin{equation}\label{A2F1KP}
\lA 
\nabla\varphi\cdot \widetilde{v}_{1}\rA_{L^{d^{\star}}}\le C_{\varphi}R,
\quad\text{with  } d^{\star}\defe 2d/(d+2).
\end{equation}
Let us assume \eqref{A2F1KP} for a moment and continue the proof. 

The constraint $\curl v_{2}=0$ implies that 
$v_{2}=\nabla{\Psi}$, for some ${\Psi}$ satisfying
\begin{equation*}
\cn (e^{\varphi}\nabla {\Psi})=-e^{\varphi}\nabla\varphi
\cdot\widetilde{v_{1}}.
\end{equation*}
This allows us to estimate $\nabla{\Psi}$ by a duality argument. We 
denote by $\scal{\cdot}{\cdot}$ the scalar product in $L^{2}$ and 
write
\begin{equation*}
\scal{e^{\varphi}\nabla{\Psi}}{\nabla{\Psi}}
=\scal{e^{\varphi}\nabla\varphi\cdot\widetilde{v}_{1}}{{\Psi}}.
\end{equation*}
Denote by $\overline{d}$ the conjugate exponent of $d^{\star}$, 
$\overline{d}\defe d^{\star}/(d^{\star}-1)=2d/(d-2)$. 
The Holder's inequality yields
\begin{equation*}
\scal{e^{\varphi}\nabla{\Psi}}{\nabla{\Psi}}
\le \lA e^{\varphi}\nabla\varphi\cdot\widetilde{v}_{1}\rA_{L^{d^{\star}}}\lA
{\Psi}\rA_{L^{\overline{d}}}.
\end{equation*}
The first factor is estimated 
by means of the claim~\eqref{A2F1KP}. 
In light of the 
Sobolev's inequality $\lA {\Psi}\rA_{L^{\overline{d}}}
\les\lA \nabla{\Psi}\rA_{L^{2}}$, we obtain
\begin{equation*}
\scal{e^{\varphi}\nabla{\Psi}}{\nabla{\Psi}}
\le C_{\varphi}R \lA \nabla{\Psi}\rA_{L^{2}}.
\end{equation*}
By using the elementary estimate 
$\lA \nabla{\Psi}\rA_{L^{2}}^{2}\le \lA e^{-\varphi}\rA_{L^{\infty}}
\scal{e^{\varphi}\nabla{\Psi}}{\nabla{\Psi}}$, we get
\begin{equation}\label{A2F1L2B}
\lA v_{2}\rA_{L^{2}}=\lA \nabla{\Psi}\rA_{L^{2}}\le C_{\varphi}R.
\end{equation}
The end of the proof is straightforward. We write
$$
\Delta {\Psi}=e^{-\varphi}\cn(e^{\varphi}\nabla{\Psi})
-\nabla\varphi\cdot\nabla{\Psi}
=-\nabla \varphi\cdot\widetilde{v}_{1}
-\nabla\varphi\cdot\nabla{\Psi},
$$
to obtain, for all $\sigma\in [0,s-1]$,
\begin{equation*}
\lA\nabla{\Psi}\rA_{H^{\sigma+1}}\les
\lA\nabla\Psi\rA_{L^{2}}+\lA\Delta{\Psi}\rA_{H^{\sigma}}\les
\lA\nabla\varphi\cdot\widetilde{v}_{1}\rA_{H^{\sigma}}
+(1+\lA\varphi\rA_{H^{s}})\lA\nabla{\Psi}\rA_{H^{\sigma}}.
\end{equation*}
To estimate the first term on the right-hand side, we verify that 
the analysis establishing~\eqref{inter} also yields
\begin{equation*}
\lA \nabla \varphi\cdot\widetilde{v}_{1}\rA_{H^{s-1}}
\les\lA\varphi\rA_{H^{s}}
\lA\nabla\widetilde{v}_{1}\rA_{H^{s-1}}\le C_{\varphi}R,
\end{equation*}
hence, by induction on $\sigma$,
\begin{equation*}
\lA\nabla{\Psi}\rA_{H^{s}}
\le C_{\varphi}R+ C_{\varphi}\lA\nabla{\Psi}\rA_{L^{2}}.
\end{equation*}
Exactly as above, one has
\begin{align*}
\lA\nabla{\Psi}\rA_{H^{s+1}}
&\les\lA\nabla\Psi\rA_{L^{2}}+\lA\Delta{\Psi}\rA_{H^{s}}\les 
\lA\nabla\varphi\cdot\widetilde{v}_{1}\rA_{H^{s}}
+(1+\lA\varphi\rA_{H^{s+1}})\lA\nabla{\Psi}\rA_{H^{s}}\\
\lA \nabla \varphi\cdot\widetilde{v}_{1}\rA_{H^{s}}
&\les\lA\varphi\rA_{H^{s+1}}
\lA\nabla\widetilde{v}_{1}\rA_{H^{s-1}}\le \lA\varphi\rA_{H^{s+1}}C_{\varphi}R.
\end{align*}
As a consequence, we end up with
\begin{equation*}
\lA\nabla{\Psi}\rA_{H^{s+1}}
\le \lA\varphi\rA_{H^{s+1}} (C_{\varphi}R+C_{\varphi}\lA\nabla{\Psi}\rA_{L^{2}}).
\end{equation*}
Therefore, the $L^{2}$ estimate \eqref{A2F1L2B} implies that
$$
\lA v_{2}\rA_{H^{s+1}}
=\lA\nabla{\Psi}\rA_{H^{s+1}}\le \lA\varphi\rA_{H^{s+1}}C_{\varphi}R.
$$
By combining this estimate with \eqref{A2F1PR1}, 
we find that
$$
\lA\nabla v_{1}\rA_{H^{s}}\le \lA\varphi\rA_{H^{s+1}}C_{\varphi}R.
$$ 
From the definition of $v_{1}$ and~\eqref{F.1}, we obtain the desired bound \eqref{A2F1}. 

We now have to establish the claim~\eqref{A2F1KP}. 

With $\overline{d}\defe 2d/(d-2)$ as above, 
the Sobolev's inequality and \eqref{F.1} imply that 
\begin{equation*}
\lA \widetilde{v}\rA_{L^{\overline{d}}}\les
\lA\nabla\widetilde{v}\rA_{L^{2}}\le C_{\varphi}R.
\end{equation*}
On the other hand, the H\"{o}lder's inequality yields
\begin{align*}
\lA \nabla {\varphi}\cdot \widetilde{v}
\rA_{L^{\delta}}\les\lA \nabla\varphi\rA_{L^{2}}
\lA \widetilde{v}\rA_{L^{\overline{d}}},\quad
 \text{  with  }\delta\defe \frac{2\overline{d}}{2+\overline{d}}=\frac{d}{d-1}\cdot
\end{align*}
By interpolating this estimate with 
$\lA \nabla\varphi\cdot\widetilde{v}\rA_{L^{\overline{d}}}\les
\lA\nabla\varphi\rA_{L^{\infty}}\lA\widetilde{v}\rA_{L^{\overline{d}}}$, 
we obtain
$$
\forall p\in [\delta,\overline{d}],\quad 
\lA \nabla\varphi\cdot \widetilde{v}\rA_{L^{p}}\les
\lA \nabla\varphi\rA_{L^{2}\cap L^{\infty}}\lA
\widetilde{v}\rA_{L^{\overline{d}}}\le C_{\varphi}R.
$$ 

Because $\curl v_{1}=0$, one can write 
$v_{1}=\nabla\Psi_{1}$ for some function 
$\Psi_{1}$ satisfying $\Delta\Psi_{1}=-\nabla\varphi\cdot\widetilde{v}$. 
Hence, the Calderon--Zygmund inequality and the previous estimate imply that 
$$
\lA \nabla v_{1}\rA_{L^{\delta}}=
\lA \nabla^{2}\Psi_{1}\rA_{L^{\delta}}\les 
\lA \Delta\Psi_{1}\rA_{L^{\delta}}\le C_{\varphi}R.
$$
Therefore, the Sobolev's inequality yields
$$
\lA v_{1}\rA_{L^{D }}\le C_{\varphi}R,\quad
\text{  with  }
D =\frac{\delta d}{d-\delta}=\frac{d}{d-2}\virgp
$$
hence, exactly as above, the H\"{o}lder's inequality gives
\begin{equation}\label{A2F1E3}
\forall p\in [\underline{d},\overline{d}],\quad \lA
e^{\varphi}\nabla\varphi\cdot \widetilde{v}_{1}\rA_{L^{p}} \le C_{\varphi}R,
\quad\text{with  }\underline{d}=\frac{2D}{2+D}=\frac{2d}{3d-4}\cdot
\end{equation}
The key estimate \eqref{A2F1KP} is now a consequence of the previous one. 
Indeed, the estimate \eqref{A2F1E3} applies with $p=d^{\star}\defe 2d/(d+2)$ since
$$
\forall d\ge 3,\qquad \underline{d}=\frac{2d}{3d-4}\le \frac{2d}{d+2}\le \frac{2d}{d-2}=\overline{d}.
$$
This completes the proof of \eqref{A2F1}.
\end{proof}

For later references, we will need the following version of~\eqref{A2F1}.
\begin{corollary}\label{coro:A2F}
Let $d= 1$ or $d\ge 3$ and $\xN\ni s>d/2$. 
There exists a function $\Cr$ such that, 
for all $\varphi\in H^{s+1}(\xR^{d})$ 
and all vector field $v\in H^{s+1}(\xR^{d})$, 
\begin{align}
\lA \nabla v\rA_{H^{s}}&\le 
\Cr(\lA \varphi\rA_{H^{s+1}})
\bigl(\lA\cn v\rA_{H^{s}}+\lA \curl(e^{\varphi}v)\rA_{H^{s}}\bigr).
\label{A2F2}
\end{align} 
\end{corollary}
\begin{proof}[Proof]
The case $d=1$ is obvious. If $d\ge 3$, Proposition~\ref{prop:A2F} 
(applied with $(\varphi,v)$ replaced
with $(-\varphi,e^{\varphi} v)$) yields
$$
\lA \nabla (e^{\varphi} v)\rA_{H^{s}}
\le \Cr(\lA \varphi\rA_{H^{s+1}})\bigl(\lA \cn  v\rA_{H^{s}}
+\lA \curl (e^{\varphi}v)\rA_{H^{s}}\bigr).
$$
Hence, to prove \eqref{A2F2} we need only prove that 
\begin{equation}\label{A2F2TD}
\lA \nabla v\rA_{H^{s}}\le 
\Cr(\lA\varphi\rA_{H^{s+1}})\lA \nabla (e^{\varphi} v)\rA_{H^{s}}.
\end{equation}
To do that we write 
$\partial_{i} v=e^{-\varphi}\partial_{i}(e^{\varphi}v)
-(e^{-\varphi}\partial_{i}\varphi)(e^{\varphi} v)$. 
The usual product rule~\eqref{PRS} implies that the $H^{s}$ norm of 
the first term is estimated by the right-hand side of~\eqref{A2F2TD}. 
Moving to the second term, we use the product
rule \eqref{inter} to obtain
$
\lA (e^{-\varphi}\partial_{i}\varphi) (e^{\varphi} v)\rA_{H^{s}}
\les (1+\lA\varphi\rA_{H^{s}})\lA \partial_{i} \varphi\rA_{H^{s}}
\lA \nabla( e^{\varphi} v)\rA_{H^{s-1}}.
$
This proves the desired bound \eqref{A2F2TD}.
\end{proof}
\begin{remark}\label{R2D}
The fact that Theorem~\ref{theo:main} 
precludes the case $d=2$ is a consequence of the fact that we do not know 
if~\eqref{A2F2} holds for $d=2$.
\end{remark}

\section{Uniform stability}\label{section:proof}
In this section, we prove Theorem~\ref{theo:main}. 
We follow closely the approach given in~\cite{TA2}: 
we recall the scheme of the analysis and 
indicate the points at which the argument must be adapted. 

Hereafter, we use the notations
\begin{gather*}
\pa\defn (\eps,\mu,\kappa)\in A\defn (0,1]\times [0,1]\times
[0,1],\qquad \nu\defn\sqrt{\mu+\kappa},\\
\norm{u}_{H^{\sigma+1}_{\alpha}}\defn\norm{u}_{\nh{\sigma}}
+\alpha\norm{u}_{\nh{\sigma+1}}\qquad (\alpha\ge 0,~\sigma\in\xR).
\end{gather*}
\subsection*{Step 1: a refined statement}
We first give our main result a refined form where the solutions satisfy 
the same estimates as the initial data do. 
Also, to prove estimates independent of $\mu$ and $\kappa$, an
important point is to seek the solutions in spaces which take into account an extra 
damping effect for the penalized terms.
\begin{definition}\label{defi:decoupling}
Let $T >0$, $\pa\defe(\eps,\mu,\kappa)\in [0,1]^{3}$ 
and set $\nu\defe\sqrt{\mu+\kappa}$.
The space $\Xr_{\pa}^{s}(T)$ consists of these $( p, v,\theta)\in C^{0}([0,T];H^{s}(\xR^{d}))$ 
such that
$$
\nu (p,v,\theta)\in C^{0}([0,T];H^{s+1}(\xR^{d})),\quad (\mu v,\kappa\theta)\in
L^{2}(0,T;H^{s+2}(\xR^{d})).
$$
The space $\Xr_{\pa}^{s}(T)$ is given the norm
\begin{align*}
\norm{( p, v,\theta)}_{\Xr_{\pa}^{s}(T)}&\defn
\norm{(\nabla p,\nabla v)}_{L^{\infty}_{T}(H^{s-1})} 
+ \norm{(\theta,\eps p,\eps v)}_{L^{\infty}_{T}(\nhsc{s+1})}\\[0.5ex]
&\quad +\sqrt{\mu}\norm{\nabla v}_{L^{2}_{T}(H^{s+1}_{\eps\nu})}+
\sqrt{\kappa} \norm{\nabla\theta}_{L^{2}_{T}(\nhsc{s+1})}\\
&\quad+ \sqrt{\mu+\kappa}\norm{\nabla p}_{L^{2}_{T}(H^{s})}
+\sqrt{\kappa}\lA\cn v\rA_{L^{2}_{T}(H^{s})},
\end{align*}
with $\lA\cdot\rA_{L^{p}_{T}(X)}$ denoting 
the norm in $L^{p}(0,T;X)$.
\end{definition}
The hybrid norm $\lA \cdot\rA_{H^{s+1}_{\eps\nu}}$ was already used by
Danchin in~\cite{Dan1}. 

For the study of nonlinear problems, it is important to 
relax the assumption that $ Q\in C^{\infty}_{0}$.
\begin{definition}
The space $F^{s}$ consists of these function 
$ Q$ such that, 
for all $\xN\ni m\le s$,
$
\partial_{t}^{m} Q\in C^{0}_{b}(\xR;H^{s+1-2m}(\xR^{d}))$, 
where $C^{0}_{b}$ stands for $C^{0}\cap L^{\infty}$.
\end{definition}

Given a normed space~$X$, we set $B(X;M)\defe \{ x\in X \,:\, \lA x\rA\le M\}$.
\begin{theorem}\label{theo:main2}
Assume that $d=1$ or $d\ge 3$ and let $\xN\ni s>1+d/2$. Given 
$M_{0}>0$ and $Q\in F^{s}$, there exist $T>0$ and $M>0$ 
such that, for all $a=(\eps,\mu,\kappa)\in\PA$ and all initial data 
$(p_{0},v_{0},\theta_{0})\in H^{s+1}(\xR^{d})$ 
satisfying 
\begin{equation}\label{I2RS}
\norm{(\nabla p_{0},\nabla v_{0})}_{H^{s-1}} 
+ \norm{(\theta_{0},\eps p_{0},\eps v_{0})}_{H^{s+1}}\le M_{0},
\end{equation}
the Cauchy problem for~\eqref{system:NSint} has a unique solution 
$(p,v,\theta)\in B(\Xr_{\pa}^{s}(T);M)$.
\end{theorem}
This theorem implies 
Theorem~\ref{theo:main}. 

\begin{remark}
A close inspection of the proof indicates 
that Theorem~\ref{theo:main2} remains valid with 
\eqref{I2RS} replaced by
$$
\norm{( p_{0}, v_{0},\theta_{0})}_{\Xr_{\pa}^{s}(0)}\defn
\norm{(\nabla p_{0},\nabla v_{0})}_{H^{s-1}} 
+ \norm{(\theta_{0},\eps p_{0},\eps v_{0})}_{H^{s+1}_{\nu}}\le M_{0}.
$$
\end{remark}
\subsection*{Step 2: local well posedness}
We explain here how to reduce matters to proving uniform bounds.
To do so, our first task is to establish the local well posedness of the Cauchy
problem for fixed $\pa=(\eps,\mu,\kappa)\in \PA$.
\begin{lemma}\label{prop:fixeda}
Let $d\ge 1$, $s>1+d/2$ and $a\in\PA$. For all 
initial data $U_{0}\defe(p_{0},v_{0},\theta_{0})\in H^{s}(\xR^{d})$, 
there exists a positive time $T$ such
that the Cauchy problem for~\eqref{system:NSint} has a unique solution
$U\defe(p,v,\theta)\in C^{0}([0,T];H^{s})$ such that $U(0)=U_{0}$. 
Moreover, the interval $[0,T^{\star})$, with 
$T^{\star}<+\infty$, is a maximal interval of $H^{s}$ existence 
if and only if $\limsup_{t\rightarrow T^{\star}}
\norm{U(t)}_{W^{1,\infty}(\xR^{d})} =+\infty$.
\end{lemma}

Lemma~\ref{prop:fixeda} is a special case of Proposition~\ref{prop:WP} 
established below. 

As in \cite{TA2,MS1}, 
on account of the previous local existence result for 
fixed $\pa\in\PA$, 
Theorem~\ref{theo:main} 
is a consequence of the following uniform estimates:

\begin{proposition} \label{proposition1} 
Let $d=1$ or $d\ge 3$, $\xN\ni s > 1+d/2$ and $M_0> 0$. There exist a constant $C_0$ 
and a non-negative function $C(\cdot )$ such that, for all 
$T \in (0,1]$ and all $\pa \in \PA$, if 
$(p,v,\theta) \in C^{\infty}([0,T];H^{\infty}(\xR^{d}))$ is a solution of~\eqref{system:NSint} with
initial data satisfying~\eqref{I2RS},
then the norm $\Omega_{\pa}(T)\defn \lA U\rA_{\Xr^{s}_{\pa}(T)}$
satisfies 
\begin{equation}\label{desi}
\Omega_{\pa}(T) \leqslant C_0 \exp \bigl((\sqrt{T}+\varepsilon)C( \Omega_{\pa}(T))\bigr).
\end{equation}
\end{proposition} 
\begin{nota}\label{nota:Notations} 
From now on, we consider an integer $s>1+d/2$, a fixed time $0<T\le 1$, a fixed 
triple of parameters $\pa=(\eps,\mu,\kappa)\in\PA$, a 
bound $M_{0}$, a fixed smooth solution $U=(p,v,\theta)
\in C^{\infty}([0,T];H^{\infty}(\xR^{d}))$ of~\eqref{system:NSint} with
initial data satisfying~\eqref{I2RS} and we set 
\begin{equation*}
\Omega\defn \lA U\rA_{\Xr^{s}_{\pa}(T)}.
\end{equation*}
With these notations, Proposition~\ref{proposition1} 
can be formulated concisely as follows: if $d\neq 2$, 
there exist constants $C_{0}$ depending only on $M_{0}$ and 
$C$ depending only on $\Omega$ such that 
$$
\Omega\le C_{0}e^{(\sqrt{T}+\eps)C}.
$$
Hereafter, we use the notations 
$\phi\defn(\theta,\eps p)$ and $\nu\defn\sqrt{\mu+\kappa}$.
\end{nota}

\begin{nota}
For later application to the nonlinear case when $Q=F(Y)$ for some unknown function $Y$, 
we also give precise estimates in terms of norms of~$Q$. 
For our purposes, the requested norm is the following:
\begin{equation}\label{defi:Qr}
\Sigma
\defn \sum_{0\le m\le s}
\lA (\id-(\eps\nu)^{2}\Delta)^{-m/2}\bigl(\eps(\partial_{t}+v\cdot\nabla)\bigr)^{m}
Q\rA_{L^{\infty}(0,T;H^{s+1-m}_{\nu})}.
\end{equation}
\end{nota}
\begin{remark}\label{Banach}
To use nonlinear estimates, it is easier to work in Banach algebras. 
If $d\ge 3$, Proposition~\ref{prop:Linftycontrol} shows that we can supplement the $\Xr^{s}_{\pa}$
estimates with $L^{\infty}$ estimates for the velocity: it suffices to prove~\eqref{desi} with $C( \Omega_{\pa}(T))$ replaced
by $C( \Omega_{\pa}^{+}(T))$ where $\Omega_{\pa}^{+}(T)\defn
\Omega_{\pa}(T)+\lA v\rA_{L^{\infty}((0,T)\times \xR^{d})}$. 
Similarly, if $d\ge 3$, all the estimates involving the source term
$Q$ remain valid with $\Sigma$ replaced by
\begin{equation*}
\sum_{0\le m\le s}
\lA (\id-(\eps\nu)^{2}\Delta)^{-m/2}(\eps\partial_{t})^{m}
Q\rA_{L^{\infty}(0,T;H^{s+1-m}_{\nu})}.
\end{equation*}
\end{remark}

\subsection*{Step 3: An energy estimate for linearized equations}
A key step in the analysis is to estimate the solution 
$(\widetilde{p},\widetilde{v},\widetilde{\theta})$ 
of linearized equations. 
As will be apparent in a moment, a notable fact is that we can 
see unsigned large terms $\eps^{-1}f^{\eps}(t,x)$ in the equations for $p$ and $v$ 
as source terms provided that: 1) they do not convey fast oscillations in time: 
$\partial_{t}f^{\eps}=O(1)$; 
2) it does not implies a loss of derivatives. 
To be more precise: in the nonlinear estimates, 
we will see the term $\eps^{-1}\chi_{1}(\phi)Q$ as a source term. Similarly, 
we can see terms of the form $\eps^{-1}F(\eps p,\theta,\sqrt{\kappa}\nabla\theta)$
 as source terms. As a result, 
it is sufficient to consider the following linearized system:
\begin{equation}\label{system:NSi}
\left\{
\begin{aligned}
&g_{1}(\phi)(\partial_{t}\widetilde{p}+ v\cdot\nabla\widetilde{p})
+\frac{1}{\eps}\cn\widetilde{v}-\frac{\kappa}{\eps}\cn(k_{1}(\phi)\nabla\widetilde{\theta})=F_{1},\\
&g_{2}(\phi)(\partial_{t}\widetilde{v}+ v\cdot\nabla\widetilde{v})
+\frac{1}{\eps}\nabla\widetilde{p}-\mu 
B_{2}(\phi,\partial_{x})\widetilde{v}=F_{2},\\
&g_{3}(\phi)(\partial_{t}\widetilde{\theta}
+ v\cdot\nabla\widetilde{\theta})
+G(\phi,\nabla\phi)\cdot\widetilde{v}+\cn\widetilde{v} -\kappa\chi_{3}(\phi)
\cn(k(\phi)\nabla\widetilde{\theta})=F_{3},
\end{aligned}
\right.
\end{equation}
where the unknown $(\widetilde{p},\widetilde{v},\widetilde{\theta})$ is 
a smooth function of $(t,x)\in [0,T]\times\xR^{d}$.

The following result establishes estimates on
\begin{equation}\label{defiL2norm}
\begin{split}
\lVert(\widetilde{p},\widetilde{v},\widetilde{\theta})\rVert_{\pa,T}
&\defn
\lA(\widetilde{p},\widetilde{v})\rA_{L^{\infty}_{T}(H^{1}_{\eps\nu})}
+\lVert\widetilde{\theta}\rVert_{L^{\infty}_{T}(H^{1}_{\nu})}\\
&\quad+\sqrt{\pe}\lVert\nabla\widetilde{\theta}\rVert_{L^{2}_{T}(H^{1}_{\nu})}
+\sqrt{\mu}\lA\nabla\widetilde{v}\rA_{L^{2}_{T}(H^{1}_{\eps\nu})}\\
&\quad+\sqrt{\mu+\kappa}\lA\nabla\widetilde{p}\rA_{L^{2}_{T}(L^{2})}
+\sqrt{\kappa}\lA\cn\widetilde{v}\rA_{L^{2}_{T}(L^{2})},
\end{split}
\end{equation}
in terms of the norm~$\lVert(\widetilde{p},\widetilde{v},\widetilde{\theta})\rVert_{\pa,0}\defn
\norm{(\widetilde{p},\widetilde{v})(0)}_{H^{1}_{\eps\nu}}
+\norm{\widetilde{\theta}(0)}_{H^{1}_{\nu}}$ of the data.

\begin{theorem}\label{theo:L2} 
Let $d\ge 1$ and assume that $G$, $k_{1}$ and $k_{3}$ are 
$C^{\infty}$ functions such that, for all $(\vartheta,\wp)\in\xR^{2}$, 
$0<k_{1}(\vartheta,\wp)<\chi_{3}(\vartheta,\wp)k(\vartheta)$. 
Set 
\begin{equation*}
R_{0}\defn\lA {\phi(0)}\rA_{H^{s-1}},\quad 
R\defn\sup_{t\in[0,T]} \lA{(\phi,\partial_{t}\phi+v\cdot\nabla\phi,
\nabla\phi,\nu\nabla^{2}\phi,\nabla v)(t)}
\rA_{H^{s-1}}.
\end{equation*}
There exist constants $C_{0}$ depending only on $R_{0}$ and 
$C$ depending only on $R$ such that,
\begin{equation*}
\norm{(\widetilde{p},\widetilde{v},\widetilde{\theta})}_{\pa,T}\le
C_{0}e^{TC}\norm{(\widetilde{p}_{0},\widetilde{v}_{0},\widetilde{\theta}_{0})}_{\pa,0}
+C\int_{0}^{T}\norm{(F_{1},F_{2})}_{H^1_{\eps\nu}}+\norm{F_{3}}_{H^{1}_{\nu}}\,dt.
\end{equation*}
\end{theorem}
In~\cite{TA2} we established the previous theorem with $R_{0}$ and $R$ 
replaced by
$$
R_{0}'=\lA\phi(0)\rA_{L^{\infty}},\quad
R'=\sup_{t\in[0,T]}
\lA(\phi,\partial_{t}\phi,v,\nabla\phi,\nu\nabla^{2}\phi,\nabla v)(t)\rA
_{L^{\infty}}.
$$
To prove the above variant, we need only check two facts. Firstly, 
in the proof of 
Theorem~$4.3$ in~\cite{TA2}, the terms $\partial_{t}\phi$ and $v$ always come together within terms 
involving the convective derivative $\partial_{t}\phi+v\cdot\nabla\phi$. 

Secondly, we have to verify that the $L^{\infty}_{t,x}$ norms of the coefficients ($g_{i}(\phi)$,...) 
are estimated by constants of the form $C_{0}e^{TC}$. In~\cite{TA2} 
we used the estimate
$$
\sup_{t\in[0,T]}\lA F(\phi(t))\rA_{L^{\infty}}
\le \lA F(\phi(0))\rA_{L^{\infty}}+T\sup_{t\in [0,T]}
\lA \partial_{t}F(\phi(t))\rA_{L^{\infty}}\le C_{0}'+TC',
$$
for some constants depending only on $R_{0}'$ and $R'$. 
Here, based on an usual estimate for hyperbolic equations, we can prove 
a similar bound:
\begin{lemma}\label{lemma:estihyp}
Let $F\in C^{\infty}(\xR^{2})$ be such that $F(0)=0$. 
There exist constants 
$C_{0}$ depending only on $R_{0}$ and $C$ depending only on $R$ such that, 
for all $t\in[0,T]$, $\lA F(\phi(t))\rA_{H^{s-1}}\le C_{0}e^{TC}$.
\end{lemma}
\begin{proof}
Since $s-1>d/2$, the Moser's estimates \eqref{PTE} and \eqref{PTE2} 
imply that there exists a function 
$\Cr$ depending only on the function $F$ such that
\begin{align*}
&\lA (\partial_{t}+v\cdot\nabla)F(\phi)\rA_{H^{s-1}}\\
&\quad\le \bigl(1+\lA F'(\phi)-F'(0)\rA_{H^{s-1}}\bigr)
\lA (\partial_{t}+v\cdot\nabla)\phi(t)\rA_{H^{s-1}},\\
&\quad\le  \Cr(\lA (\phi,\partial_{t}\phi+v\cdot\nabla\phi)\rA_{H^{s-1}})\le \Cr(R),
\end{align*}
and $\lA F(\phi(0))\rA_{H^{s-1}}\le \Cr(\lA \phi(0)\rA_{H^{s-1}})$.

Hence, the desired estimate follows from the following 
estimate: there exists a constant $V$ depending only 
on $\lA\nabla v\rA_{L^{\infty}_{T}H^{s-1}}$ such that
\begin{equation*}
\sup_{t\in[0,T]}\lA F(\phi(t))\rA_{H^{s-1}}
\le \lA F(\phi(0))\rA_{H^{s-1}}+TV\sup_{t\in [0,T]}
\! \!\lA (\partial_{t}+v\cdot\nabla)F(\phi(t))\rA_{H^{s-1}}.
\end{equation*}
To prove this result we set 
$\widetilde{u}\defn \partial_{x}^{\alpha}F(\phi)$ where $\alpha\in\xN^{d}$ is such that 
$\la\alpha\ra\le s-1$. Then $\widetilde{u}$ 
solves
\begin{equation*}
\partial_{t}u+v\cdot\nabla u=f\defn 
\partial_{x}^{\alpha} \bigl((\partial_{t}+v\cdot\nabla)F(\phi)\bigr)
+[v,\partial_{x}^{\alpha}]\cdot\nabla F(\phi).
\end{equation*}
Since $s-1>d/2$, the product rule~\eqref{PRS} implies that
\begin{align*}
\lVert
[v,\partial_{x}^{\alpha}]\cdot\nabla F(\phi)\rVert_{L^{2}}
&\les \sum_{\beta+\gamma=\alpha,~\beta>0}
\lVert \partial_{x}^{\beta}v \partial_{x}^{\gamma}\nabla F(\phi)\rVert_{L^{2}}\\
&\les \sum_{\beta+\gamma=\alpha,~\beta>0}
\lVert \partial_{x}^{\beta}v\rVert_{H^{s-1-(\la\beta\ra-1)}}
\lA \partial_{x}^{\gamma}\nabla F\rA_{H^{s-1-(\la\gamma\ra+1)}},
\end{align*}
hence, 
$
\lA f\rA_{L^{2}}\les 
\lA(\partial_{t}+v\cdot\nabla)F(\phi)\rA_{H^{s-1}}
+\lA \nabla v\rA_{H^{s-1}}\lA F(\phi)\rA_{H^{s-1}}$.
 
We next use an integration by parts argument yielding
$$
\frac{d}{dt}\lA\widetilde{u}\rA_{L^{2}}^{2}\le 
(1+\lA \cn v\rA_{L^{\infty}})\lA\widetilde{u}\rA_{L^{2}}^{2}
+\lA f\rA_{L^{2}}^{2}.
$$
The Gronwall's Lemma concludes the proof.
\end{proof}

\subsection*{Step 4: High frequency estimates}
On the technical side, the estimate of the derivatives is divided
into four steps. Most of the work concerns the separation of estimates
into high and low frequency components, where the division occurs at
frequencies of order of the inverse 
of $\eps \nu$ where $\nu\defe\sqrt{\mu+\kappa}$. 

We begin by estimating the high frequency component
\begin{equation*}
\Omega^{\HF}\defn \lA (\id-J_{\eps\nu})U\rA_{\Xr^{s}_{\pa}(T)},
\end{equation*}
where $\{ J_{h}\,\arrowvert\, h\in [0,1]\}$ is a Friedrichs
mollifiers: $J_{h}=\jmath(h D_{x})$ is the Fourier multiplier 
with symbol $\jmath(h \xi)$ where 
$\jmath$ is a~$C^{\infty}$ function 
of~$\xi\in\mathbb{R}^d$, satisfying
\begin{equation*}
0\le \jmath \le 1, \quad \jmath(\xi) = 1 \text{ for } |\xi| \le 1, 
\quad \jmath(\xi)=0 \text{ for } |\xi| \ge 2,\quad \jmath(\xi)
=\jmath\left(-\xi\right).
\end{equation*}

\begin{proposition}\label{PRHF} 
Let $d\ge 1$. There exist constants 
$C_0$ depending 
only on $M_{0}$ and $C$ depending only on $\Omega$, 
such that 
\begin{equation}\label{resultat3}
\Omega^{\HF}
\leqslant C_{0}e^{\sqrt{T}C}+\sqrt{T}C\lA Q\rA_{L^{\infty}_{T}(H^{s+1}_{\nu})}.
\end{equation}
\end{proposition} 
\begin{proof}[Proof]
Introduce $ P \defn(\id -J_{\eps\nu})\Fi{}{s}$ and 
$\widetilde{U}\defn (Pp,Pv,P\theta)$. 
Then, $\widetilde{U}$ satisfies 
System~\eqref{system:NSi} with
$$
k_{1}(\phi)\defn \chi_{1}(\phi)k(\theta),\qquad 
G(\phi,\nabla\phi)\defn g_{3}(\phi)\nabla\theta,
$$
and 
$F=(F_{1},F_{2},F_{3})^{T}\defn f_{\HF}+f_{Q}+f_{\chi}$, 
where 
\begin{equation*}
f_{Q}\defn
\begin{pmatrix} \eps^{-1}P \bigl(\chi_{1}(\phi)Q\bigr)\\
0\\ P \bigl(\chi_{3}(\phi)Q\bigr)\end{pmatrix},
\quad 
f_{\chi}\defn 
\begin{pmatrix}
-\kappa\eps^{-1}\nabla\chi_{1}(\phi)\cdot (k(\theta)\nabla\widetilde{\theta})\\
0\\
0\end{pmatrix},
\end{equation*}
and $f_{\HF}$ is given by
\begin{alignat*}{3}
f_{1,\HF}&\defe\bigl[g_{1}(\phi), P \bigr](\partial_{t}+v\cdot\nabla)p
&&+g_{1}(\phi)\bigl[v, P \bigr]\cdot\nabla p  
&&-\frac{\kappa}{\eps}\bigl[B_{1}(\phi,\partial_{x}), P \bigr]\theta,\\
f_{2,\HF}&\defe\bigl[g_{2}(\phi), P \bigr](\partial_{t}+v\cdot\nabla)v
&&+g_{2}(\phi)\bigl[v, P \bigr]\cdot\nabla v
&&-\mu\bigl[B_{2}(\phi,\partial_{x}), P \bigr]v,\\[0.5ex]
f_{3,\HF}&\defe \bigl[g_{3}(\phi), P \bigr](\partial_{t}+v\cdot\nabla)\theta
&&+g_{3}(\phi)\bigl\{v; P \bigr\}\cdot\nabla\theta
&&-\pe\bigl[B_{3}(\phi,\partial_{x}), P \bigr]\theta,
\end{alignat*}
where $B_{i}(\phi,\partial_{x})\defe
\chi_{i}(\phi)\cn(k(\theta)\nabla\cdot)$ ($i=1,3$), $[A,B]=AB-BA$ and
\begin{equation*}
\bigl\{v; P \bigr\}\cdot\nabla\theta
\defn v \cdot\nabla P \theta +(P v)\cdot\nabla\theta- P(v\cdot\nabla\theta).
\end{equation*}

\smallskip
\noindent{\emph{Estimate for $f_{\HF}$}}. We use the following analogue of
Lemma~$5.3$ in~\cite{TA2}: there exists a constant $K=K(d,s)$ such that
\begin{align*}
\norm{\bigl[ f, P \bigr]u}_{H^{1}_{\eps\nu}}
&\le \eps\nu K\norm{\nabla f}_{L^{\infty}}\norm{u}_{\nh{s}}
+\eps\nu K\lA \nabla f\rA_{H^{s}}\lA u\rA_{L^{\infty}},\\
\norm{\bigl[ f, P \bigr]u}_{H^{1}_{\nu}}
&\les \nu K\lA \nabla f\rA_{L^{\infty}}\lA u\rA_{\nh{s}}
+\nu K\lA\nabla f\rA_{H^{s}}\lA u\rA_{L^{\infty}}.
\end{align*}
The fact that the right-hand side only involves $\nabla f$ follows from 
the most simple of all the sharp commutator estimates established in~\cite{David}: 
for all $s>1+d/2$ and all Fourier multiplier 
$A(D_{x})\in \Op S^{s}_{1,0}$, there exists a constant $K$ 
such that, for all 
$f\in H^{s}(\xR^{d})$ and all $u\in H^{s}(\xR^{d})$, 
\begin{equation}
\lA [ f,A(D_{x})]u\rA_{L^{2}}
\le K\lA \nabla f\rA_{L^{\infty}}\lA u\rA_{H^{s-1}}+K\lA \nabla f\rA_{H^{s-1}}
\lA u\rA_{L^{\infty}}.
\end{equation}

As in \cite{TA2}, 
from this and the usual nonlinear estimates~\eqref{PTE} and \eqref{PTE2}, 
it can be verified that there exists a generic function $\Cr$ (depending only on parameters 
that are considered fixed) such that,
\begin{align*}
\lA f_{1,\HF}\rA_{H^{1}_{\eps\nu}}
&\le \Cr\bigl( \lA (\theta,\eps p,\eps v)\rA_{H^{s+1}_{\nu}}\bigr)
\{1+\lA \eps(\partial_{t}+v\cdot\nabla)p\rA_{H^{s}_{\nu}}+\kappa\lA \theta\rA_{H^{s+2}}\},\\
\lA f_{2,\HF}\rA_{H^{1}_{\eps\nu}}
&\le \Cr\bigl( \lA (\theta,\eps p,\eps v)\rA_{H^{s+1}_{\nu}}\bigr)
\{1+\lA \eps(\partial_{t}+v\cdot\nabla )v\rA_{H^{s}_{\nu}}+\mu\lA \eps v\rA_{H^{s+2}}\},\\
\lA f_{3,\HF}\rA_{H^{1}_{\nu}}
&\le \Cr\bigl( \lA (\theta,\eps p,\eps v)\rA_{H^{s+1}_{\nu}}\bigr)
\{1+\lA (\partial_{t}+v\cdot\nabla)\theta\rA_{H^{s}_{\nu}}+\kappa\lA\theta\rA_{H^{s+2}}\}\cdot
\end{align*}
Set $\psi\defe(\theta,\eps p,\eps v)$. The key point is that
\begin{equation}\label{defi:psi}
\begin{split}
&\lA (\partial_{t}+v\cdot\nabla)\psi\rA_{H^{s}_{\nu}}\\
&~~
\le \Cr(\lA \psi\rA_{H^{s+1}_{\nu}})
\bigl\{1
+\lA(\nu\nabla p,\nu\cn v,\eps\mu\nabla^{2}v,
\kappa\nabla^{2}\theta)\rA_{H^{s}}
+\lA Q\rA_{H^{s}_{\nu}}\bigr\}\cdot
\end{split}
\end{equation}
This estimate differs from the one that appears in Lemma $5.14$ in~\cite{TA2} in that 
the right-hand side does not involve $v$ itself but only its derivatives. Yet, 
as the reader can verify, the same proof applies since we do 
not estimate $\partial_{t}\psi$ but instead $\partial_{t}\psi+v\cdot\nabla\psi$.

\smallskip
\noindent\emph{Estimate for $f_{Q}$ and $f_{\chi}$}. 
By using the elementary estimate
$$
\lA (\id-J_{\eps\nu})u\rA_{H^{\sigma+1}_{\eps\nu}}\les \eps\nu\lA u\rA_{H^{\sigma+1}},
$$
we find that
\begin{align*}
\frac{1}{\eps}\lA P (\chi_{1}(\phi)Q)\rA_{H^{1}_{\eps\nu}}
+\lA P  (\chi_{3}(\phi)Q)\rA_{H^{1}_{\nu}} 
\le \lA  \chi_{1}(\phi)Q\rA_{H^{s+1}_{\nu}}+\lA  \chi_{3}(\phi)Q\rA_{H^{s+1}_{\nu}}.
\end{align*}
The tame estimates~\eqref{PTE} and \eqref{PTE2} 
(see also Lemma~$5.5$ and $5.6$ in~\cite{TA2}) imply
\begin{align*}
\lA\chi_{i}(\phi)Q\rA_{H^{s+1}_{\nu}}
\les 
(1+\lA \chi_{i}(\phi)-\chi_{i}(0)\rA_{H^{s+1}_{\nu}})\lA Q\rA_{H^{s+1}_{\nu}}
\les \Cr(\lA \phi\rA_{H^{s+1}_{\nu}})\lA Q\rA_{H^{s+1}_{\nu}}
\end{align*}
so that 
$\norm{f_{1,Q}}_{L^{\infty}_{T}(H^1_{\eps\nu})}+\norm{f_{3,Q}}_{L^{\infty}_{T}(H^{1}_{\nu})}\le
C\lA Q\rA_{L^{\infty}_{T}(H^{s+1}_{\nu})}$. 
The technique for estimating $f_{\chi}$ is similar; we find that 
$\lA f_{1,\chi}\rA_{L^{\infty}_{T}(H^1_{\eps\nu})}\le C$.

\smallskip
By definition of $\lA\cdot\rA_{\Xr_{\pa}^{s}(T)}$, the previous estimates imply 
that there exists a constant $C$ depending only on $\Omega$ such that 
\begin{align*}
\int_{0}^{T}\norm{(F_{1},F_{2})}_{H^1_{\eps\nu}}+\norm{F_{3}}_{H^{1}_{\nu}}\,dt 
&\le\sqrt{T}\Bigl(\int_{0}^{T}\norm{(F_{1},F_{2})}_{H^1_{\eps\nu}}^{2}+\norm{F_{3}}_{H^{1}_{\nu}}^{2}\,dt\Bigr)^{1/2} \\
&\le \sqrt{T}C+\sqrt{T}C\lA Q\rA_{L^{\infty}_{T}(H^{s+1}_{\nu})}.
\end{align*}

From here we can parallel the rest of the argument of Section~$5$ in~\cite{TA2}, to prove that 
$\lA ( P p, P v, P \theta)\rA_{\pa,T}\le C_{0}\exp
(\sqrt{T}C)+\sqrt{T}C\lA Q\rA_{L^{\infty}_{T}(H^{s+1}_{\nu})}$
where the norm $\lA\cdot\rA_{\pa,T}$ 
is as defined in~\eqref{defiL2norm}. 
Since $\Omega^{\HF}\les \lA ( P p, P v, P \theta)\rA_{\pa,T}$, this completes the proof.
\end{proof}

\subsection*{Step 5: Low frequency estimates}

The following step 
is to estimate the low frequency part of the
fast components:
\begin{align*}
\Omega^{\LF}&\defn 
\lA \cn J_{\eps\nu}v\rA_{L^{\infty}_{T}(H^{s-1})}
+\nu\lA \cn J_{\eps\nu}v\rA_{L^{2}_{T}(H^{s})}\\
&\quad +\lA \nabla J_{\eps\nu} p\rA_{L^{\infty}_{T}(H^{s-1})}
+\nu\lA \nabla J_{\eps\nu}p\rA_{L^{2}_{T}(H^{s})}.
\end{align*}

\begin{proposition}\label{PRLF}
Let $d\ge 1$. There exist constants 
$C_0$ depending 
only on~$M_{0}$, 
$C$ depending only on $\Omega$ 
and $C'$ depending only on $\Omega+\Sigma$, 
such that 
\begin{equation}\label{LFE1}
\Omega^{\LF}
\leqslant C_{0}e^{(\sqrt{T}+\eps)C}+\sqrt{T}C'.
\end{equation}
\end{proposition}

By contrast with the high frequency regime, the estimate~\eqref{LFE1} 
cannot be obtained 
from the $L^{2}$ estimates by an elementary argument using 
differentiation of the equations (see~\cite{MS1,Schochet}). 
To overcome this problem, we first give estimates 
for the time derivatives, and next we use 
the special structure of the equations to estimate the spatial derivatives. 

For the case of greatest physical 
interest ($d=3$), the proof 
given in~\cite{TA2} applies with only minor changes. 
Indeed, as alluded to in Remark~\ref{Banach}, 
it suffices to check that all the estimates 
involving~$\lA v\rA_{H^{s}}$ remain valid with 
$\lA v\rA_{H^{s}}$ replaced by $\lA v\rA_{L^{\infty}}+\lA \nabla 
v\rA_{H^{s-1}}$. Yet, if $d\le 2$, because of the lack of $L^{2}$ estimates 
for the velocity, we cannot 
use the time derivatives. For this problem, we use an idea 
introduced by Secchi in~\cite{Secchi}. Namely, 
we replace~$\partial_{t}$ by the convective derivative~
$$
D_{v}\defn\partial_{t}+v\cdot\nabla.
$$
For the reader convenience, we 
indicate how to adapt the three main calculus inequalities in~\cite{TA2} 
when $\partial_{t}$ is replaced by~$D_{v}$.

First, to localize in the low frequency region we use 
the following commutator estimate. The think of interest is the gain
of an extra factor~$\eps$.
\begin{lemma}\label{bignorms:commutator}
Given $s>1+d/2$, there exists a constant $K$ such that for all
$\varepsilon\in [0,1]$, all $\nu\in [0,2]$, all $T>0$, all $\ik\in\mathbb{N}$ 
such that $1\le\ik\le s$ and all 
$f,u$ and $v$ in $C^{\infty}([0,T];H^{\infty}(\mathbb{D}))$,
\begin{align*}
&\norm{\bigl[f,J_{\varepsilon\nu}(\eps D_{v})^{\ik}]u}_{H^{s-\ik+1}_{\varepsilon\nu}}\\
&\quad\le
K\varepsilon\Bigl\{\norm{f}_{H^{s}}+\sum_{\ell=0}^{\ik-1}
\norm{\Fi{\varepsilon\nu}{-\ell}(\varepsilon D_{v})^{\ell}D_{v}f}_{H^{s-1-\ell}}\Bigr\}\\
&\quad\quad\times
\Bigl\{
\norm{\Fi{\varepsilon\nu}{-m}(\varepsilon D_{v})^{m}u}_{\nhsc{s-\ik}}+\sum_{\ell=0}^{\ik-1}
\norm{\Fi{\varepsilon\nu}{-\ell}(\varepsilon D_{v})^{\ell}u }_{\nh{s-1-\ell}}\Bigr\}\virgp
\end{align*}
where $\Fi{\eps\nu}{\sigma}\defn (\id-(\eps\nu)^{2}\Delta)^{\sigma/2}$.
\end{lemma}

To apply the previous lemma, we need estimates of the coefficients 
$f$ and~$D_{v}f$. Since, for System~\eqref{system:NSint}, the coefficients are functions of the slow
variable $(\theta,\eps p,\eps v)$, the main estimates are the following.

\begin{lemma}\label{prop:notrouble}
Let $s>1+d/2$ be an integer. There exists a function~$\Cr(\cdot)$ such
that, for all $\pa=(\varepsilon,\mu,\kappa)\in A$, all $T >0$ 
and all smooth solution $(\vari,\vard,\variii)\in
C^{\infty}([0,T];H^{\infty}(\mathbb{D}))$ of~\eqref{system:NSint}, 
if $\nu \in [(\mu+\kappa)/2,2]$ then 
the function $\Psi$ defined by
\begin{equation*}
\Psi\defn \bigl(\psi,D_{v}\psi,\nabla\psi\bigr)
\quad\mbox{where}\quad \psi\defn (\variii,\varepsilon\vari,\varepsilon\vard),
\end{equation*}
satisfies 
\begin{gather}
\sum_{0\le \ell\le s}
\norm{\Fi{\varepsilon\nu}{-\ell}(\varepsilon D_{v})^{\ell}\Psi}_{\nh{s-\ell-1}}
\le \Cr \bigl(\norm{\Psi}_{H^{s-1}}+\Sigma\bigr),\label{induction:TER1}\\
\sum_{0\le \ell\le s}\norm{\Fi{\varepsilon\nu}{-\ell}(\varepsilon D_{v})^{\ell}\Psi}_{\nhsc{s-\ell}}
\le \Cr \bigl(\norm{\Psi}_{H^{s-1}}+\Sigma\bigr)\norm{\Psi}_{\nhsc{s}},\label{induction:TER2}
\end{gather}
where $\Sigma$ is as defined in~$\eqref{defi:Qr}$.
\end{lemma}

Once this is granted, we are in position to estimate 
the commutator of the equations~\eqref{system:NSint} 
and $\mathcal{P}\defn J_{\eps\nu}(\eps D_{v})^{s}$: 
\begin{alignat*}{3}
f_{1,\LF}&\defe\bigl[g_{1}(\phi), \mathcal{P} \bigr]D_{v}p
&&+g_{1}(\phi)\bigl[v, \mathcal{P} \bigr]\cdot\nabla p  
&&-\frac{\kappa}{\eps}\bigl[B_{1}(\phi,\partial_{x}), \mathcal{P} \bigr]\theta,\\
f_{2,\LF}&\defe\bigl[g_{2}(\phi), \mathcal{P} \bigr]D_{v}v
&&+g_{2}(\phi)\bigl[v, \mathcal{P} \bigr]\cdot\nabla v
&&-\mu\bigl[B_{2}(\phi,\partial_{x}), \mathcal{P} \bigr]v,\\[0.5ex]
f_{3,\LF}&\defe \bigl[g_{3}(\phi), \mathcal{P} \bigr]D_{v}\theta
&&+g_{1}(\phi)\bigl[v, \mathcal{P} \bigr]\cdot\nabla \theta  
&&-\pe\bigl[B_{3}(\phi,\partial_{x}), \mathcal{P} \bigr]\theta.
\end{alignat*}
It is found that 
$$
\lA f_{1,\LF}\rA_{H^{1}_{\eps\nu}}
+\lA f_{1,\LF}\rA_{H^{1}_{\eps\nu}}
+\lA f_{1,\LF}\rA_{H^{1}_{\nu}}
\le (1+\lA\Psi\rA_{H^{s}_{\nu}})\Cr(\lA \Psi\rA_{H^{s-1}}+\Sigma).
$$
Note that $\Psi$ is estimated by means of~\eqref{defi:psi}.

As in the high frequency regime, we have to estimate 
source terms of the form~$\eps^{-1}\mathcal{P} F(\Psi,Q)$. 
The fact that these large source terms cause 
no difficulty comes from the fact that 
$\eps^{-1}J_{\eps\nu}(\eps D_{v})^{s}F(\Psi,Q)=
J_{\eps\nu}(\eps D_{v})^{s-1}D_{v}F(\Psi,Q)$ 
together with $D_{v}F(\Psi,Q)=O(1)$ 
(the norm~$\Sigma$ introduced in~\eqref{defi:Qr} is the requested 
norm to give this statement a precise meaning). 

With these results in hands, one can estimate 
$J_{\eps\nu}(\eps D_{v})^{s}(p,v,\theta)$ by means of 
Theorem~\ref{theo:L2}. 
Next, we give estimate for $\cn J_{\eps\nu}v$ and $\nabla J_{\eps\nu}p$ 
from the estimate of $J_{\eps\nu}(\eps D_{v})^{s}(p,v,\theta)$ by
means of the following induction argument:
\begin{lemma}
Set $\norm{u}_{\mathcal{K}^{\sigma}_{\nu}(T)}\defn
\norm{u}_{L^{\infty}(0,T;H^{\sigma-1})}
+\nu\norm{u}_{L^{2}(0,T;H^{\sigma})}$.

Let $\widetilde{\var}\defn (\tvari,\tvard,\tvariii)$ solve
\begin{equation}\label{Xsystem:NS5}
\left\{
\begin{aligned}
&g_{1}(\phi)(\partial_{t}\tvari+v\cdot\nabla\tvari)
+\eps^{-1}\cn \tvard
-\pe\eps^{-1}\chi_{1}(\phi)\cn(k(\theta)\nabla\tvariii)=f_{1},\\
&g_{2}(\phi)(\partial_{t}\tvard+v\cdot\nabla\tvard)
+\eps^{-1}\nabla\tvari-\r B_{2}(\phi,\partial_{x})\tvard=f_{2},\\
&g_{3}(\phi)(\partial_{t}\tvariii+v\cdot\nabla\tvariii)
+\cn\tvard -\pe \chi_{3}(\phi)\cn(k(\theta)\nabla\tvariii)=f_{3}.
\end{aligned}
\right.
\end{equation}
If support of the Fourier transform of $\tvar$ is included in the 
ball~$\{\la\xi\ra\le 2/\eps\cpe\}$, then there exist constant $C_{0}$
depending only on $M_{0}$ and $C$ depending only on $\Omega$ 
such that, for all $\sigma \in [1,s]$,
\begin{equation}\label{Xesti:LFinduction}
\begin{split}
&\norm{\nabla\tvari}_{\Kr_{\cpe}^{\sigma}(T)} 
+\norm{\cn\tvard}_{\Kr_{\cpe}^{\sigma}(T)}\\[0.5ex]
&\quad\le \widetilde{C}\norm{(\eps D_{v})\tvari}_{\Kr_{\cpe}^{\sigma}(T)} 
+\widetilde{C}\norm{(\eps D_{v})\cn\tvard}_{\Kr_{\cpe}^{\sigma-1}(T)}\\[0.5ex]
&\quad\quad +\widetilde{C}\norm{\nabla\tvari}_{L^{\infty}_{T}(L^{2})}+
\widetilde{C} \norm{\tvariii(0)}_{\nhsc{\sigma+1}}
+\eps C\norm{\r\tvard}_{\Kr_{\cpe}^{\sigma+1}(T)}\\[0.5ex]
&\quad\quad + \eps C\norm{(f_{1},f_{2})}_{\Kr_{\cpe}^{\sigma}(T)}
+\cpe \widetilde{C}\norm{f_{3}}_{L^{2}_{T}(H^{\sigma})},
\end{split}
\end{equation}
where $\widetilde{C}\defn C_{0}e^{(\sqrt{T}+\eps)C}$.
\end{lemma}

\subsection*{Step 6: estimates for the slow components.}
To complete the proof of~\eqref{desi}, it remains to estimate $\curl v$
and $\theta$. Yet, this is not straightforward. Following 
M\'etivier and Schochet~\cite{MS1}, we begin by estimating 
$\curl( \gamma v)$ for some appropriate positive weight $\gamma=\Gamma(\theta,\eps p)$.

\begin{lemma}\label{LSC} 
Let $d\ge 1$. There exist constants 
$C_0$ depending 
only on $M_{0}$ and $C$ depending only on $\Omega$, and there 
exists a function $\Gamma\in C^{\infty}(\xR^{2})$ such that, with 
$\gamma\defe\Gamma(\theta,\eps p)$, there holds
\begin{equation*}
\lA \curl (\gamma v)\rA_{L^{\infty}_{T}(H^{s-1})}
+\sqrt{\mu}\lA \curl (\gamma v)\rA_{L^{2}_{T}(H^{s})}\le 
C_{0}e^{\sqrt{T}C}+\sqrt{T}C\lA Q\rA_{L^{\infty}_{T}(H^{s+1}_{\nu})}.
\end{equation*}
\end{lemma}

\begin{lemma}\label{LSTHETA}
Let $d\ge 1$. There exist constants 
$C_0$ depending 
only on $M_{0}$ and $C$ depending only on $\Omega$, such that
\begin{equation*}
\lA J_{\eps\nu}\theta\rA_{L^{\infty}_{T}(H^{s+1}_{\nu})}
+\sqrt{\kappa}\lA J_{\eps\nu}\theta\rA_{L^{2}_{T}(H^{s+2}_{\nu})}
\le C_{0}e^{\sqrt{T}C}+\sqrt{T}C\lA Q\rA_{L^{\infty}_{T}(H^{s+1}_{\nu})}.
\end{equation*}
\end{lemma} 
The proofs of Lemma~\ref{LSC} and \ref{LSTHETA} follow 
from a close inspection of the proofs of
Lemma~$6.25$ and $6.26$ in~\cite{TA2}. 
We just mention that this is where we use 
the function $\Di$ of Assumption~(H\ref{assus}) in~\S\ref{assu:assus} ($\gamma$ is related to the fluid entropy).

\begin{lemma}\label{lemma:curl}
Assume $d\ge 3$. There exist constants 
$C_0$ depending 
only on $M_{0}$ and $C$ depending only on $\Omega$ 
such that, with $\gamma_{0}\defe\Gamma(\theta_{0},\eps p_{0})$ where 
$\Gamma$ is as above, there holds
\begin{equation*}
\lA \curl (\gamma_{0}v)\rA_{L^{\infty}_{T}(H^{s-1})}
+\sqrt{\mu}\lA\curl(\gamma_{0}v)\rA_{L^{2}_{T}(H^{s})}\le C_{0}e^{\sqrt{T}C}+\sqrt{T}C\lA Q\rA_{L^{\infty}_{T}(H^{s+1}_{\nu})}.
\end{equation*}
\end{lemma}
\begin{proof}
Set $\widetilde{\gamma}\defn\gamma-\gamma_{0}$. 
By~Lemma~\ref{LSC}, all we need to prove is that 
\begin{equation}\label{desiWG}
\lA \curl (\widetilde{\gamma}v)\rA_{L^{\infty}_{T}(H^{s-1})}
+\sqrt{\mu}\lA \curl (\widetilde{\gamma}v)\rA_{L^{2}_{T}(H^{s})}\le 
\sqrt{T}C+\sqrt{T}C\lA Q\rA_{L^{\infty}_{T}(H^{s+1}_{\nu})}.
\end{equation}
To do so, we claim that $\widetilde{\gamma}$ is small for small times:
\begin{equation}\label{WTG}
\lA \widetilde{\gamma}\rA_{L^{\infty}_{T}(H^{s})}+\nu\lA\widetilde{\gamma}\rA_{L^{2}_{T}(H^{s+1})}
\le \sqrt{T}C+\sqrt{T}C\lA Q\rA_{L^{\infty}_{T}(H^{s+1}_{\nu})}.
\end{equation}
Let us assume this and continue the proof. 

We have to estimate 
$\curl (\widetilde{\gamma}v)=\widetilde{\gamma}\curl v +(\nabla\widetilde{\gamma})\times v$. 
By combining the Cauchy-Schwarz estimate with 
the usual product rule~\eqref{PRS} and the product rule~\eqref{inter}, we find that
\begin{align*}
\lA\widetilde{\gamma}\curl v\rA_{L^{\infty}_{T}(H^{s-1})}&
\le \lA\widetilde{\gamma}\rA_{L^{\infty}_{T}(H^{s-1})}\lA \nabla v\rA_{L^{\infty}_{T}(H^{s-1})},\\
\sqrt{\mu}\lA\widetilde{\gamma}\curl v\rA_{L^{2}_{T}(H^{s})}&
\le \lA\widetilde{\gamma}\rA_{L^{\infty}_{T}(H^{s})}\lA \sqrt{\mu}\nabla v\rA_{L^{2}_{T}(H^{s})},\\
\lA\nabla\widetilde{\gamma}\times v\rA_{L^{\infty}_{T}(H^{s-1})}&
\le \lA\widetilde{\gamma}\rA_{L^{\infty}_{T}(H^{s})}\lA \nabla v\rA_{L^{\infty}_{T}(H^{s-1})},\\
\sqrt{\mu}\lA\nabla\widetilde{\gamma}\times v\rA_{L^{2}_{T}(H^{s})}
&\le \lA\sqrt{\mu}\widetilde{\gamma}\rA_{L^{2}_{T}(H^{s+1})}\lA \nabla v\rA_{
L^{\infty}_{T}(H^{s-1})}.
\end{align*}
The claim~\eqref{WTG} then yields the desired bound~\eqref{desiWG}.

We now have to prove the claim~\eqref{WTG}. We first note that
\begin{align*} 
\nu\lA\widetilde{\gamma}\rA_{L^{2}_{T}(H^{s+1})}&\le
\nu\sqrt{T}\lA\widetilde{\gamma}\rA_{L^{\infty}_{T}(H^{s+1})}\\
&\le \nu\sqrt{T}\Cr(\lA (\theta,\eps p)\rA_{L^{\infty}_{T}(L^{\infty}_{x})})
(1+\lA (\theta,\eps p)\rA_{L^{\infty}_{T}(H^{s+1})})\\
&\le \sqrt{T}\Cr(\lA (\theta,\eps p)\rA_{L^{\infty}_{T}(H^{s+1}_{\nu})})\le \sqrt{T}C.
\end{align*}
To prove the second half of~\eqref{WTG}, we verify that, directly from the definitions, 
$\widetilde{\gamma}$ satisfies an equation of the form 
$\partial_{t}\widetilde{\gamma}+v\cdot\nabla\widetilde{\gamma}=f$ with $f$ 
bounded in $L^{2}(0,T;H^{s}(\xR^{d}))$ by a constant depending only
on~$\Omega+\lA Q\rA_{L^{\infty}_{T}(H^{s+1}_{\nu})}$. Then, 
we apply the above mentioned estimate for hyperbolic
equations: 
\begin{equation}
\norm{\widetilde{\gamma}}_{L^{\infty}_{T}(H^{s})}\lesssim
e^{TV}\lA\widetilde{\gamma}(0)\rA_{H^{s}}
+\int_{0}^{T}e^{(T-t)V}\lA f\rA_{H^{s}}\,dt,
\end{equation}
where $V\defe K\int_{0}^{T}\norm{\nabla v}_{H^{s-1}}\,dt$ with $K=K(s,d)$. Since $\widetilde{\gamma}(0)=0$, 
by applying the Cauchy-Schwarz inequality, it is found 
that the $L^{\infty}_{T}(H^{s})$ norm of~$\widetilde{\gamma}$ is estimated by 
$\sqrt{T}e^{TV}\lA f\rA_{L^{2}_{T}(H^{s})}$, 
thereby obtaining the claim.
\end{proof}

\subsection*{Step 7: closed set of estimates}
To complete the proof of Proposition~\ref{proposition1}, 
it remains to check that we have proved a closed set of estimates. 

The obvious estimate 
$\lA u\rA_{H^{\sigma}}\le \lA J_{\eps\nu}u\rA_{H^{\sigma}}+
\lA (\id-J_{\eps\nu})u\rA_{H^{\sigma}}$ implies that
$$
\lA (\nabla p,\cn v)\rA_{L^{\infty}_{T}(H^{s-1})} 
+ \sqrt{\mu+\kappa}\lA(\nabla p,\cn v)\rA_{L^{2}_{T}(H^{s})}
\les \Omega_{\LF}+\Omega_{\HF},
$$
and, similarly, $\lA{\theta}\rA_{L^{\infty}_{T}(H^{s+1}_{\nu})}
+\sqrt{\kappa} \lA{\nabla\theta}\rA_{L^{2}_{T}(\nhsc{s+1})}$ is
estimated by
\begin{equation*}
\lA J_{\eps\nu}{\theta}\rA_{L^{\infty}_{T}(H^{s+1}_{\nu})}
+\sqrt{\kappa} \lA J_{\eps\nu}{\nabla\theta}\rA_{L^{2}_{T}(\nhsc{s+1})}
+\Omega_{\HF}.
\end{equation*}
The estimate $\lA \eps u\rA_{H^{\sigma+1}_{\nu}}
\les \lA\eps u\rA_{L^{2}}+\lA \nabla u\rA_{H^{\sigma-1}}+\lA
(\id-J_{\eps\nu})u\rA_{H^{\sigma+1}_{\nu}}$ yields
\begin{equation}\label{esti:nablav}
\begin{split}
&\lA (\eps p,\eps v)\rA_{L^{\infty}_{T}(H^{s+1}_{\nu})}+
\sqrt{\mu}\lA \nabla v\rA_{L^{2}_{T}(H^{s+1}_{\eps\nu})}\les \lA (\eps
p,\eps v)\rA_{L^{\infty}_{T}(L^{2})} \\
&\qquad\qquad\qquad+\lA (\nabla p,\nabla v)\rA_{L^{\infty}_{T}(H^{s-1})}+
\sqrt{\mu}\lA \nabla v\rA_{L^{2}_{T}(H^{s})}+\Omega_{\HF}.
\end{split}
\end{equation}
On the other hand, Corollary~\ref{coro:A2F} implies that, if $d\neq
2$, there exists a constant $C_{0}$ depending only on $M_{0}$ such that
\begin{align*}
&\lA\nabla v\rA_{L^{\infty}_{T}(H^{s-1})}+\sqrt{\mu}\lA \nabla v\rA_{L^{2}_{T}(H^{s})}\\
&\qquad\le C_{0}\lA (\cn v,\curl(\gamma_{0}v))\rA_{L^{\infty}_{T}(H^{s-1})}+C_{0}\sqrt{\mu}
\lA (\cn v,\curl(\gamma_{0}v))\rA_{L^{2}_{T}(H^{s})}.
\end{align*}
By using the estimate~\eqref{defi:psi}, one can verify that the term $\lA (\eps p,\eps v)\rA_{L^{\infty}_{T}(L^{2})}$
(in the left-hand side of~\eqref{esti:nablav}) can be
estimated as in the proof of Lemma~\ref{lemma:estihyp}. Therefore, 
according to Propositions~\ref{PRHF}--\ref{PRLF} and 
Lemma~\ref{LSTHETA}--\ref{lemma:curl}, 
we have proved that, if $d\neq 2$, then 
$\Omega\le\widetilde{C}$ where $\widetilde{C}=C_{0}e^{(\sqrt{T}+\eps)C}+\sqrt{T}C'$ for some
constants $C_{0}$, $C$ and $C'$ depending only on $M_{0}$, $\Omega$
and $\Omega+\Sigma$, respectively. 

This concludes the proof of Proposition~\ref{proposition1} and 
hence Theorem~\ref{theo:main2}.

\section{Uniform estimates in the Sobolev spaces}\label{section:Sobolev}
With regards to the low Mach number limit problem, 
we mention that the convergence results
\footnote{These results are strongly based on a Theorem 
of M\'etivier and Schochet~\cite{MS1} about the decay 
to zero of the local energy for a class of wave operators 
with variable coefficients.}  proved in~\cite{TA2} apply for general systems (not only for 
perfect gases). To avoid repetition, we only mention that one can rigorously 
justify the low Mach number limit for general initial data provided that one can prove that the 
solutions are uniformly bounded in 
Sobolev spaces (see Proposition $8.2$ in~\cite{TA2}). The problem presents itself: 
Theorem~\ref{theo:main} 
only gives uniform estimates for the derivatives of~$p$ and $v$. 
In this section, we give uniform bounds in Sobolev norms.

\begin{theorem}\label{theo:GP}
Let $d\ge 1$ and $\xN\ni s>1+d/2$. 
Assume that $Q=0$. Also, assume that either 
$\chi_{1}=\chi_{1}(\vartheta,\wp)$ is independent of 
$\vartheta$ or that $d\ge 3$. Then, for all $M_{0}>0$, there exists $T>0$ and $M>0$ 
such that, for all $a=(\eps,\mu,\kappa)\in\PA$ and all initial data 
$(p_{0},v_{0},\theta_{0})\in H^{s+1}(\xR^{d})$ satisfying
$$
\norm{(p_{0},v_{0},\theta_{0})}_{\nh{s+1}}\le M_{0},
$$
the Cauchy problem for~\eqref{system:NSint} has a unique 
classical solution $(p,v,\theta)$ in
$C^{0}([0,T];H^{s+1}(\xR^{d}))$  such that
$$
\sup_{t\in [0,T]}\lA (p(t),v(t),\theta(t))\rA_{H^{s}}\le M.
$$
\end{theorem}
The first half of this result is proved in \cite{TA2}. 
Indeed, the assumption that $\chi_{1}(\vartheta,\wp)$ does not depend on $\vartheta$ 
is satisfied by perfect gases. So we concentrate on the second 
half ($d\ge 3$). In view of Theorem~\ref{theo:main2}, 
it remains only to prove 
{\em a posteriori\/} uniform~$L^{2}$ estimates. 
More precisely, the proof of Theorem~\ref{theo:GP} reduces to
establishing the following result.
\begin{lemma}\label{EUL2}
Let $d\ge 3$. Consider a family of solutions 
$(p^{\pa},v^{\pa},\theta^{\pa})$ of~$\eqref{system:NSint}$ (for some source terms $Q^{\pa}$) 
uniformly bounded 
in the sense of the conclusion of Theorem~$\ref{theo:main2}$:
\begin{equation}\label{defi:EL2WS}
\sup_{\pa\in\PA}
\lA ( p^{\pa},v^{\pa},\theta^{\pa})\rA_{\Xr_{\pa}^{s}(T)}<+\infty,
\end{equation}
for some $s>1+d/2$ and fixed $T>0$. 
Assume further that the source terms $Q^{\pa}$ are uniformly bounded in 
$C^{1}([0,T];L^{1}\cap L^{2}(\xR^{d}))$ and that the initial data 
$(p^{\pa}(0),v^{\pa}(0))$ are uniformly bounded in $L^{2}(\xR^{d})$. Then 
the solutions $(p^{\pa},v^{\pa},\theta^{\pa})$ are uniformly bounded in $C^{0}([0,T];L^{2}(\xR^{d}))$.
\end{lemma}
\begin{remark}
We allow $Q^{\pa}$ for application to the combustion equations. 
To clarify matters, we note that one can replace~\eqref{defi:EL2WS} by
$$
\sup_{\pa\in\PA}\sup_{t\in [0,T]}\lA (\nabla p^{\pa}(t),
\nabla v^{\pa}(t))\rA_{H^{s}}+\lA \theta^{\pa}(t)\rA_{H^{s+1}}<+\infty,
$$
for some $s>2+d/2$.
\end{remark}
\begin{proof}
For this proof, we set
\begin{align*}
R\defn \sup_{\pa\in\PA}
\Bigl\{\lA ( p^{\pa},v^{\pa},\theta^{\pa})\rA_{\Xr_{\pa}^{s}(T)}+
\lA (p^{\pa}(0),v^{\pa}(0))\rA_{L^{2}}
+\lA Q^{\pa}\rA_{C^{1}([0,T];L^{1}\cap L^{2})}\Bigr\},
\end{align*}
and we denote by $C(R)$ various constants depending only on~$R$.

The strategy of the proof consists of transforming the system~\eqref{system:NSint} 
so as to obtain $L^{2}$ estimates uniform in $\eps$ by a simple integration 
by parts argument. 

To do that we claim that there exist $U^{\pa}\in C^{1}([0,T];L^{2}(\xR^{d}))$
satisfying the following properties: 
\begin{align}
\sup_{\pa\in\PA}\lA (p^{\pa},v^{\pa})\rA_{L^{\infty}_{T}(L^{2})}&\le 
\sup_{\pa\in\PA}\lA U^{\pa}\rA_{L^{\infty}_{T}(L^{2})}+C(R),\label{EUL2:H1}\\
\sup_{\pa\in\PA}\lA U^{\pa}(0)\rA_{L^{2}}&\le C(R),\label{EUL2:H15}
\end{align}
and $U^{\pa}$ solves a system having the form 
\begin{equation}\label{EUL2:2}
E^{\pa}\partial_{t}U^{\pa}+\eps^{-1}S(\partial_{x})U^{\pa}=F^{\pa},
\end{equation}
where $S(\partial_{x})$ is 
skew-symmetric, the matrices $E^{\pa}=E^{\pa}(t,x)$ are positive
definite and one has the uniform bounds
\begin{equation}\label{EUL2:H2}
\sup_{\pa\in\PA}
\lA \partial_{t}E^{\pa}\rA_{L^{\infty}([0,T]\times\xR^{d})}
+\lA (E^{\pa})^{-1}\rA_{L^{\infty}([0,T]\times\xR^{d})}^{-1}
+\lA F^{\pa}\rA_{L^{1}_{T}(L^{2})}\le C(R).
\end{equation}

\smallbreak
Before we prove the claim, let us prove that it implies Lemma~\ref{EUL2}. 
To see this, we combine two basic ingredients: 
\begin{align*}
\frac{d}{dt}\scal{E^{\pa} U^{\pa}}{U^{\pa}}
&=-\eps^{-1}\scal{S(\partial_{x})U^{\pa}}{U^{\pa}}+\scal{F^{\pa}}{U^{\pa}}
+\scal{(\partial_{t} E^{\pa})U^{\pa}}{U^{\pa}}\\
&\le \lA F^{\pa}\rA_{L^{2}}^{2}+C(R)\lA U^{\pa}\rA_{L^{2}}^{2},
\end{align*}
and 
$\lA U^{\pa}\rA_{L^{2}}^{2}
\le\lA (E^{\pa})^{-1}\rA_{L^{\infty}}^{-1}\scal{E^{\pa} U^{\pa}}{U^{\pa}}$. Hence, 
by~\eqref{EUL2:H15} and~\eqref{EUL2:H2}, the Gronwall's Lemma implies that $\lA U^{\pa}\rA_{L^{\infty}_{T}(L^{2})}\le C(R)$. 
The estimate \eqref{EUL2:H1} thus implies the desired result.

To prove the claim, we set $U^{\pa}\defn (p^{\pa},v^{\pa}-V^{\pa})^{T}$ where
$$
V^{\pa}\defn
\kappa\chi_{1}(\phi^{\pa})k(\theta^{\pa})\nabla\theta^{\pa}
+\nabla\Delta^{-1}\bigl(-\kappa\nabla\chi_{1}(\phi^{\pa})\cdot k(\theta^{\pa})\nabla\theta^{\pa}
+\chi_{1}(\phi^{\pa})Q^{\pa}\bigr).
$$
The fact that $V^{\pa}$ is well defined follows
from~Proposition~\ref{prop:RD}. We next verify that $U^{\pa}$ satisfies~\eqref{EUL2:2} with 
$$
E^{\pa}=\begin{pmatrix}g_{1}(\phi^{\pa}) & 0 \\
0 & g_{2}(\phi^{\pa})\end{pmatrix}, \quad
S(\partial_{x})=\begin{pmatrix} 0 & \cn \\ \nabla & 0
\end{pmatrix},
$$
$$
F^{\pa}=\begin{pmatrix}
-g_{1}(\phi^{\pa})v^{\pa}\cdot\nabla p^{\pa}\\
-g_{1}(\phi^{\pa})v^{\pa}\cdot\nabla v^{\pa}+\mu B_{2}(\phi^{\pa},\partial_{x})v^{\pa}
-g_{2}(\phi^{\pa})\partial_{t}V^{\pa}\end{pmatrix}.
$$

By~\eqref{PRS},~\eqref{PTE2} and~\eqref{restriction:dimge3}, 
to prove that the bounds~\eqref{EUL2:H1} and \eqref{EUL2:H2} hold, 
it suffices to prove 
that $\lA \partial_{t}\phi^{\pa}\rA_{H^{s-1}}\le C(R)$. Yet, this is nothing new. 
Indeed, we first observe that, directly from the equations, 
$$
\lA\partial_{t}\phi^{\pa}+v^{\pa}\cdot\nabla\phi^{\pa} \rA_{H^{s-1}}\le C(R).
$$
On the other hand, the product rule~\eqref{inter} implies that
$\lA v^{\pa}\cdot\nabla\phi^{\pa}\rA_{H^{s-1}}$ is estimated by 
$\lA \nabla v^{\pa}\rA_{H^{s-1}}\lA \phi^{\pa}\rA_{H^{s}}\le C(R)$. 
This completes the proof.
\end{proof}

\begin{remark}For our purposes, one of the main differences between 
$\xR^{3}$ and~$\xR$ is the following. 
For all $f\in C^{\infty}_{0}(\xR^{3})$, 
Proposition~\ref{propRD1} implies that there exists a vector field $u\in 
H^{\infty}(\xR^{3})$ such that $\cn u=f$. In sharp contrast, 
the mean value of the divergence of a smooth vector field $u\in H^{\infty}(\xR)$ 
is zero. This implies 
that Lemma~\ref{EUL2} is false with $d=1$.
\end{remark}

The following result contains an analysis of the easy case where initially $\theta_{0}=O(\eps)$. 
This regime is interesting for the 
incompressible limit (see \cite{BN}).
\begin{proposition}\label{prop:WP}
Let $d\ge 1$ and $\xR\ni s>1+d/2$. For all $M_{0}>0$, there exists $T>0$ and $M>0$ 
such that for all $a\in\PA$ and all initial data 
$(p_{0},v_{0},\theta_{0})\in H^{s}(\xR^{d})$ satisfying
\begin{equation}\label{PIA1}
\lA (p_{0},v_{0})\rA_{\nh{s}}+\eps^{-1}\lA \theta_{0}\rA_{H^{s}}\le M_{0},
\end{equation}
the Cauchy problem for~\eqref{system:NSint} has a unique classical 
solution $(p,v,\theta)$ in
$C^{0}([0,T];H^{s}(\xR^{d}))$  such that
\begin{equation}\label{PWPUB}
\sup_{t\in [0,T]}\lA (p(t),v(t))\rA_{H^{s}}+\eps^{-1}\lA \theta (t)\rA_{H^{s}}\le M.
\end{equation}
\end{proposition}
\begin{proof}
The proof of this result is based on the change of unknown 
$( p, v,\theta)\mapsto
(\DII(\theta,\eps p), v,\theta)$ where 
$\DII$ is as given by Assumption~(H\ref{assus}) in~$\S$\ref{assu:assus}. 
By setting $\rho\defe \DII(\theta, \eps p)$ 
it is found that $(p,v,\theta)$ satisfies \eqref{system:NSint} 
if and only if 
\begin{equation}\label{Psystem:MHPa}
\left\{
\begin{aligned}
&\chi_{3}(\partial_{t}\rho+\vard\cdot\nabla\rho) + (\chi_{3}-\chi_{1})\cn\vard=0,\\
&g_{2}(\ffp{\vard})+\varepsilon^{-2}\gamma_{1}\nabla\variii+\varepsilon^{-2}\gamma_{2}\nabla\rho
-  \mu B_{2}\vard  = 0,\\
&g_{3}(\ffp{\variii})+\cn\vard -  \kappa \chi_{3}\cn(k\nabla \variii)=0,
\end{aligned}
\right.
\end{equation}
where $\gamma_{1}\defe(\chi_{1}g_{3})/(\chi_{3}g_{1})$ and $\gamma_{2}\defe 1/g_{1}$. 
Notice that Assumption (H\ref{assus}) implies that the 
coefficients 
$g_{i}$, $\gamma_{i}$, $\chi_{3}$ and $\chi_{3}-\chi_{1}$ are positive.

The key point is that the assumption~\eqref{PIA1} allows us to symmetrize the 
equations by setting $u\defn (\tilde{\rho},v,\tilde{\theta})$, where
$$
\tilde{\rho}\defn \eps^{-1}\rho,\quad 
\tilde{\theta}\defn \eps^{-1}\theta.
$$ 
The fact that this change of unknowns is singular in $\eps$ 
causes no difficulty. Indeed, directly from the assumption \eqref{PIA1}, we have 
$\norm{\tilde{\theta}(0)}_{H^{s}}\le M_{0}$. On the other hand, 
the assumption $\DII(0,0)=0$ implies that there is a function $C_{G}$ such that 
$\lA \DII(u)\rA_{H^{\sigma}}\le C_{\DII}(\lA u\rA_{L^{\infty}})\lA u\rA_{H^{\sigma}}$ 
for all $u\in H^{\sigma}$ with 
$\sigma> d/2$. Therefore, we have
\begin{equation}\label{MDIIR}
\begin{split}
\lA \tilde{\rho}\rA_{H^{s}}&=\eps^{-1}\lA \DII(\theta,\eps p)\rA_{H^{s}}
\le \eps^{-1}C_{\DII}(\lA (\theta,\eps p)\rA_{L^{\infty}})\lA (\theta,\eps p)\rA_{H^{s}}\\
&\le C_{\DII}(\lA (\theta,\eps p)\rA_{L^{\infty}})\lVert (\tilde{\theta},p)\rVert_{H^{s}},
\end{split}
\end{equation}
hence, $\lA \tilde{\rho}(0)\rA_{H^{s}}\le C_{0}$ for 
some constant depending only on $M_{0}$.

Because $(\vartheta,\wp)\mapsto (\vartheta,\DII(\vartheta,\wp))$ is a $C^{\infty}$ diffeomorphism 
with $\DII(0,0)=0$, one can write $\eps p=P\bigl(\theta,G(\theta,\eps p)\bigr)=P(\theta,\rho)$,
for some $C^{\infty}$ function $P$ vanishing at the origin. Therefore
one can see the coefficients ($g_{i}$, $\chi_{i}$, $\gamma_{i}$...) as
functions of $(\theta,\rho)$. Hence, with
$u=(\tilde{\rho},v,\tilde{\theta})$ as above, 
one can rewrite System \eqref{Psystem:MHPa} 
under the form
\begin{equation}\label{PE1}
A_{0}(\eps u)\partial_{t}u +\sum_{1\le j\le d}
A_{j}(u,\eps u)\partial_{j}u
+
\frac{1}{\eps}\sum_{1\le j\le d}S_{j}(\eps u)\partial_{j}u
-B(\eps u,\partial_{x})u=0,
\end{equation}
where the matrices $S_{j}$, $A_{j}$ are symmetric (with 
$A_{0}$ positive definite) 
and the viscous perturbation 
$B(\eps u,\partial_{x})$ is as in \eqref{Psystem:MHPa}. 
Furthermore, one can always assume 
that the matrices $S_{j}$ have constant coefficients. 

Since the matrix $A_{0}$ multiplying the time derivative 
depends only on the unknown through $\eps u$, 
and since the initial data $u(0)$ are uniformly bounded in $H^{s}$, 
the proof of the uniform existence Theorem of~\cite{KM2} applies. By that proof, 
we conclude that the solutions of \eqref{PE1} exist and are 
uniformly bounded for a time $T$ 
independent of $\eps$. 
Once this is granted, it remains to verify that 
the solutions $(p,v,\theta)$ of System~\eqref{system:NSint} 
exist and are uniformly bounded in the sense of~\eqref{PWPUB}. 
To see this, as for $\tilde{\rho}$ in \eqref{MDIIR}, we note that
\begin{align*}
\lA p\rA_{H^{s}}&=\lA P(\theta,\rho)\rA_{H^{s}}\\
&\le \eps^{-1}C_{P}(\lA (\theta,\rho)\rA_{L^{\infty}})\lA (\theta,\rho)\rA_{H^{s}}
= C_{P}(\lA (\theta,\rho)\rA_{L^{\infty}})\lVert (\tilde{\theta},\tilde{\rho})\rVert_{H^{s}}\\
&\le C(\lVert (\tilde{\theta},\tilde{\rho})\rVert_{H^{s}}),
\end{align*}
so that 
$\lA (p,v)\rA_{H^{s}}+\eps^{-1}\lA \theta \rA_{H^{s}}\le 
C(\lA u\rA_{H^{s}})$. This completes the proof.
\end{proof}
\begin{remark}
Consider the Euler equations ($\mu=0=\kappa$ and $\eps=1$). 
By a standard re-scaling, Proposition~\ref{prop:WP} just says that the classical 
solutions with small initial data of size $\delta$ exist 
for a time of order of $1/\delta$. 
Following the approach initiated by Alinhac in~\cite{Alinhac}, several much more 
precise results have been obtained. 
In particular, 
the interested reader is referred to the recent advance of Godin~\cite{Godin_ARMA} 
(for the $3$D nonisentropic Euler equations). 
\end{remark}

\section{Spatially periodic solutions}\label{section:Torus}
In this section, we consider the case where $x$ belongs to the torus $\xT^{d}$. 
\begin{theorem}\label{theo:mainP}
Let $d\ge 1$ and $\xN\ni s>1+d/2$. 
For all source term $Q\in C^{\infty}(\xR\times\xT^{d})$ and for all $M_{0}>0$, there exist $T>0$ and $M>0$ 
such that, for all $a\in\PA$ and all initial data 
$(p_{0},v_{0},\theta_{0})\in H^{s+1}(\xT^{d})$ satisfying
$$
\lA(p_{0},v_{0})\rA_{H^{s}} 
+ \norm{(\theta_{0},\eps p_{0},\eps v_{0})}_{\nh{s+1}}\le M_{0},
$$
the Cauchy problem for~\eqref{system:NSint} has a unique classical solution 
$(p,v,\theta)$ in $C^{0}([0,T];H^{s+1}(\xT^{d}))$ such that
$$
\sup_{t\in [0,T]}\lA\nabla p(t)\rA_{H^{s-1}}+\lA  v(t)\rA_{H^{s}} 
+ \norm{(\theta(t),\eps p(t))}_{\nh{s}}\le M.
$$
\end{theorem}
The proof follows from two observations: first, 
the results proved in Steps~$1$--$6$ in section~\ref{section:proof} 
apply {\em mutatis mutandis} in the periodic case; and second, 
as proved below, 
the periodic case is easier in that one can prove uniform $L^{2}$ estimates for the 
velocity. This in turn implies that (as in~\cite{TA2,MS1}) 
one can directly prove a closed set of estimates 
by means of the estimate:
$$
\lA v\rA_{H^{s}(\xT^{d})}\le C\lA\cn v\rA_{H^{s-1}(\xT^{d})}
+C\lA \curl(\gamma v)\rA_{H^{s-1}(\xT^{d})}
+C\lA v\rA_{L^{2}(\xT^{d})},
$$
for some constant $C$ depending only on $\lA
\log\gamma\rA_{H^{s}(\xT^{d})}$ (compare with~\eqref{A2F2}). 

Let us concentrate on the main new qualitative property:
\begin{lemma}\label{EUL3}
Let $d\ge 1$. Consider a family of solutions 
$(p^{\pa},v^{\pa},\theta^{\pa})$ of~$\eqref{system:NSint}$ (for some source terms $Q^{\pa}$) 
such that
$$
\sup_{\pa\in\PA}\lA
(p^{\pa},v^{\pa},\theta^{\pa})\rA_{\Xr_{\pa}^{s}(T)}<+\infty,
$$
for some $s>1+d/2$ and fixed $T>0$. 
If $Q^{\pa}$ is uniformly bounded in 
$C^{1}([0,T];L^{2}(\xT^{d}))$ and 
$(p^{\pa}(0),v^{\pa}(0))$ is uniformly bounded in $L^{2}(\xT^{d})$, then 
$v^{\pa}$ is uniformly bounded in $C^{0}([0,T];L^{2}(\xT^{d}))$.
\end{lemma}
\begin{proof}
The main new technical ingredient is, as used by Schochet
in~\cite{SchoV}, an appropriate ansatz for the pressure. 

Again, the proof makes use of the Fourier multiplier $\nabla\Delta^{-1}$. 
Note, that $\nabla\Delta^{-1}$ is bounded from
$L^{2}_{\sharp}(\xT^{d})$ to $H^{1}(\xT^{d})$ where
$L^{2}_{\sharp}(\xT^{d})$ consists of these functions $u\in L^{2}(\xT^{d})$ such that 
$\L{u}\defn\int_{\xT^{d}}u(x)\,dx=0$. 

Set 
$$
F^{\pa}\defn \kappa\chi_{1}(\phi^{\pa})\cn(k(\theta^{\pa})\nabla\theta^{\pa})
+\chi_{1}(\phi^{\pa})Q^{\pa},
$$
and introduce the functions 
$V^{\pa}=V^{\pa}(t,x)$ and $P^{\pa}=P^{\pa}(t)$ by
$$
P^{\pa}\defn \frac{\L{F^{\pa}}}{\L{g_{1}(\phi^{\pa})}}
\quad\text{and}\quad
V^{\pa}\defn \nabla\Delta^{-1}\bigl(F^{\pa}-g_{1}(\phi^{\pa})P^{\pa}\bigr),
$$
so that
$$
F^{\pa}=g_{1}(\phi^{\pa})P^{\pa}+\cn V^{\pa}.
$$ 
This allows us to rewrite 
the first equation in~\eqref{system:NSint} as 
$$
g_{1}(\phi^{\pa})
(\partial_{t}p^{\pa}+v^{\pa}\cdot\nabla p^{\pa})
+\eps^{-1}\cn (v^{\pa}-V^{\pa})=g_{1}(\phi^{\pa})P^{\pa}.
$$ 
Therefore, by introducing 
$$
U^{\pa}\defn (q^{\pa},v^{\pa}-V^{\pa})^{T}\quad\text{with}\quad 
q^{\pa}(t,x)=p^{\pa}(t,x)-P^{\pa}(t),
$$ 
we are back in the situation of Lemma~\ref{EUL2}: $U^{\pa}$ satisfies
\begin{equation}\label{EUL2:7}
E^{\pa}(\partial_{t}U^{\pa}+v^{\pa}\cdot\nabla U^{\pa})+\eps^{-1}S(\partial_{x})U^{\pa}=F^{\pa},
\end{equation}
where $S(\partial_{x})$ is 
skew-symmetric, the matrices $E^{\pa}$ are positive definite and
\begin{equation*}
\lA (E^{\pa},\partial_{t}E^{\pa}+v^{\pa}\cdot\nabla E^{\pa})\rA_{L^{\infty}([0,T]\times\xR^{d})}
+\lA (E^{\pa})^{-1}\rA_{L^{\infty}([0,T]\times\xR^{d})}^{-1}
+\lA F^{\pa}\rA_{L^{1}_{T}(L^{2})}
\end{equation*}
is uniformly bounded. 

As before, the proof proceeds by multiplying by
$U^{\pa}$ and integrating on $[0,T]\times\xT^{d}$. We find that 
$\partial_{t}\scal{E^{\pa} U^{\pa}}{U^{\pa}}$ is given by
\begin{equation*}
\scal{\bigl((\partial_{t}+v^{\pa}\cdot\nabla)E^{\pa}\bigr)U^{\pa}}{U^{\pa}}+\scal{E^{\pa}(\cn
  v^{\pa})U^{\pa}}{U^{\pa}}
+2\scal{F^{\pa}}{U^{\pa}}\virgp
\end{equation*}
and hence conclude that 
$U^{\pa}$ is uniformly bounded in $C^{0}([0,T];L^{2}(\xT^{d}))$. 
Since $V^{\pa}$ is uniformly bounded in $C^{0}([0,T];L^{2}(\xT^{d}))$, 
this yields the desired result.
\end{proof}
\begin{remark}
In the periodic case, 
as shown by M\'etivier and Schochet~\cite{MS2,MS3} as well as Bresch, 
Desjardins, Grenier and Lin~\cite{BDGL}, 
the study of the behavior of the solutions when $\eps\rightarrow 0$ 
involved many additional phenomena. 
\end{remark}
\section{Low Mach number combustion}\label{section:Combustion}
The system \eqref{system:ANS} is relevant whenever all nuclear 
or chemical reactions are frozen, which is the case in many treatments
of fluid mechanics. By contrast, for the combustion, one has
to replace the energy evolution equation by
$$
\partial_{t}(\rho e)+\cn(\rho v e)+P\cn
 v=\kappa\cn(k\nabla{T})+F(Y),
$$
with $Y\defn (Y_{1},\ldots,Y_{L})$ where the $Y_{\ell}$'s 
denote the relative
concentrations of nuclear or chemical species. The new unknown 
$Y_{\ell}$ satisfies~:
\begin{equation}
\partial_{t}(\rho Y_{\ell})
+\cn(\rho v Y_{\ell})=\lambda\cn(D_{\ell}\nabla Y_{\ell})+\rho \omega_{\ell}(t,x),
\end{equation}
where $\omega_{\ell}$ is a given source term, $D_{\ell}>0$ and $\lambda$ measures 
the importance of diffusion processes. 

Many results have been obtained for the reactive gas equations 
(see~\cite{CHT} and the references therein). Yet, the previous
studies do not include the dimensionless numbers.
Here we consider the system:
\begin{equation}\label{system:ANSF3}
\left\{
\begin{aligned} 
&\alpha(\partial_{t}P+ v\cdot\nabla P) +\cn v
=\kappa\beta\cn(k\nabla{T})+F_{1}(Y,T,P),\\
&\rho(\partial_{t} v+ v\cdot\nabla v)+\frac{\nabla
  P}{\eps^{2}}=\mu\bigl(2\cn(\zeta D v)+\nabla(\eta\cn v)\bigr),\\
&\gamma(\partial_{t}{T}+ v\cdot\nabla{T})+\cn v=
\kappa\delta\cn(k\nabla{T})+F_{3}(Y,T,P),\\[0.5ex]
&\rho(\partial_{t} Y
+v\cdot\nabla Y)=\lambda\cn(D\nabla Y),
\end{aligned}
\right.
\end{equation}
where $\alpha$, $\beta$, $\gamma$ and $\delta$ are given 
functions of $(Y,T,P)$. 

As explained in the introduction, it is convenient to introduce
$(p,\theta,y)$ by 
$P=\underline{P}e^{\eps p}$, 
$T=\underline{T}e^{\theta}$, $Y=\underline{Y}e^{y}$, 
where $(\underline{P},\underline{T},\underline{Y})\in [0,+\infty)^{2+L}$. 
For smooth solutions, $(p,v,\theta,y)$ satisfies a system of the form:
\begin{equation}\label{system:NSint3}
\left\{
\begin{aligned}
&g_{1}(\Phi)(\partial_{t}p+v\cdot\nabla p)
+\frac{1}{\eps}\cn v =
\frac{\kappa}{\eps}\chi_{1}(\Phi)\cn(k(\theta)\nabla\theta)+\frac{1}{\eps}Q_{1}(\Phi),\\
&g_{2}(\Phi)(\partial_{t}v +v\cdot\nabla v) +\frac{1}{\eps}\nabla p 
= \mu \chi_{2}(\Phi)\bigl(\cn (\zeta(\theta) D v)
+\nabla(\eta(\theta) \cn v)\bigr),\\
&g_{3}(\Phi)(\partial_{t}\theta +v\cdot\nabla \theta)+\cn v
=\kappa \chi_{3}(\Phi)\cn(k(\theta)\nabla\theta)+Q_{3}(\Phi),\\[0.5ex]
&g_{4}(\Phi)(\partial_{t}y +v\cdot\nabla y)
=\lambda\chi_{4}(\Phi)\cn(D(\theta)\nabla y),
\end{aligned}
\right.
\end{equation}
where $\Phi\defe (y,\theta,\eps p)$.

\begin{assu}Denote by $(\mathsf{y},\vartheta,\wp)\in\xR^{N}$ the 
place holder of the unknown $(y,\theta,\eps p)$. 
Parallel to Assumption~(H\ref{assus}) in~$\S$\ref{assu:assus}, we
suppose that $\gi_{i}$ and $\bi_{i}$ ($i=1,2,3$) are $C^{\infty}$ 
positive functions of $(\mathsf{y},\vartheta,\wp)\in\xR^{N}$, $\bi_{1}<\bi_{3}$ 
and there exist two functions $\Di$ and $\DII$ such that  
$(\mathsf{y},\vartheta,\wp)\mapsto(\Di(\mathsf{y},\vartheta,\wp),\wp)$ and 
$(\mathsf{y},\vartheta,\wp)\mapsto(\mathsf{y},\vartheta,\DII(\vartheta,\wp))$ 
are $C^{\infty}$ diffeomorphisms from $\xR^{N}$ onto $\xR^{N}$, $F$ and $G$ vanish at the origin, 
and
\begin{equation*}
\gi_{1}\frac{\partial \Di}{\partial \vartheta}
=-\gi_{3}\frac{\partial \Di}{\partial \wp} >0,\qquad 
\gi_{1}\bi_{3}\frac{\partial \DII}{\partial \vartheta}
=-\gi_{3}\bi_{1}\frac{\partial \DII}{\partial \wp} <0.
\end{equation*}
Moreover, $Q_{1}$ and $Q_{3}$ are $C^{\infty}$ 
functions of $(\mathsf{y},\vartheta,\wp)$ 
vanishing at the origin.
\end{assu}

Introduce
$$
B\defn \bigl\{ \,(\eps,\mu,\kappa,\lambda)\in 
(0,1]\times [0,1]\times [0,1]\times [0,2]\,\arrowvert\,
\lambda\ge \sqrt{\mu+\kappa}\,\bigr\}\cdot
$$

\begin{definition}\label{defi:decoupling2}
Let $T >0$, $s\in\xR$, $b=(\eps,\mu,\kappa,\lambda)\in B$ and set $\pa\defn(\eps,\mu,\kappa)$. 
The space $\Zr_{b}^{s}(T)$ consists of these $( p, v,\theta,y)\in C^{0}([0,T];H^{s}(\xR^{d}))$ 
such that
$$
(p,v,\theta)\in\Xr_{\pa}^{s}(T),\quad
\nu y\in C^{0}([0,T];H^{s+1}(\xR^{d})),
\quad \lambda y\in L^{2}(0,T;H^{s+2}_{\nu}(\xR^{d})),
$$
where $\nu\defn \sqrt{\mu+\kappa}$ and $\Xr_{\pa}^{s}(T)$ is 
as defined in Definition~\ref{defi:decoupling}. 
The space $\Zr_{b}^{s}(T)$ is given the norm
\begin{align*}
\norm{(p, v,\theta,y)}_{\Zr_{b}^{s}(T)}&\defn
\norm{(p, v,\theta)}_{\Xr_{\pa}^{s}(T)}
+\lA y\rA_{L^{\infty}_{T}(H^{s+1}_{\nu})}
+\sqrt{\lambda}\lA y\rA_{L^{2}_{T}(H^{s+2}_{\nu})}.
\end{align*}
\end{definition}
Having proved estimates for the solutions of System~\eqref{system:NSint} with precised
estimates in terms of the norm~$\Sigma$ of the source term $Q$ (see~\eqref{defi:Qr}), 
we are now in position to assert that:
\begin{theorem}\label{theo:main3}
Assume that $d\neq 2$. Given 
$M_{0}>0$ and $\xN\ni s>1+d/2$, there exist $T>0$ and $M>0$, 
such that for all $b\in B$ and all initial data 
$(p_{0},v_{0},\theta_{0},y_{0})\in H^{s+1}(\xR^{d})$ satisfying
$$
\norm{(\nabla p_{0},\nabla v_{0})}_{H^{s-1}} 
+ \norm{(y_{0},\theta_{0},\eps p_{0},\eps v_{0})}_{H^{s+1}}\le M_{0},
$$
the Cauchy problem for~\eqref{system:NSint} has a unique classical solution 
$(p,v,\theta,y)$ in the ball $B(\Zr_{b}^{s}(T);M)$.
\end{theorem}
\begin{remark}
For the case of greatest physical interest ($d=3$), 
Theorem~\ref{theo:main3} has two corollaries. 
As alluded to in Section~\ref{section:Sobolev}, 
it allows us to rigorously justify, at least in the whole space case, 
the computations given by Majda in~\cite{Maj}. By the way, 
this proves the well posedness of the Cauchy problem for 
the zero Mach number combustion in the whole space 
(this was known only in the periodic case~\cite{Embid1}). 
Moreover, note that the solutions given by Theorem~\ref{theo:main3} satisfy
uniform estimates recovering in the limit~$\eps\rightarrow 0$ those
obtained by Embid for the limit system. 
Finally, we mention that the previous analysis seems to apply with 
$Q_{i}(\Phi)$ replaced by $\chi_{i}(\Phi)Q(\Phi,\nabla y,\nabla^{2}y)$
for some smooth function $Q$, yet we will not address this issue.
\end{remark}

\appendix

\section{General equations of state}\label{appendix:COV}
Recall that, in order to study the full Navier-Stokes
equations~\eqref{system:ANS}, 
we choose to work with the unknown $(P,v,{T})$. 
In order to close this system, we must relate $(\rho,e)$ to $(P,T)$ 
by means of two equations of state: $\rho=\rho(P,T)$ and $e=e(P,T)$. 
The purpose of this section is to show that 
Assumption~(H\ref{assus}) in~$\S$\ref{assu:assus} is satisfied under
general assumptions on the partial derivatives of $\rho$ and $e$
with respect to $P$ and~$T$. 

\subsection{Computation of the coefficients}\label{SA1CC}
We begin by expressing the 
coefficients $g_{i}$ and $\chi_{i}$, which appear in
\eqref{system:NSint}, in terms of 
the partial derivatives of $\rho$ and $e$
with respect to $P$ and~$T$. To do that it is convenient 
to introduce the entropy. Here is where the first identity in 
\eqref{P1FI} enters.
\begin{assu}\label{assu:general}
The functions $\rho$ and $e$ are $C^{\infty}$ functions 
of~$(P,{T})\in (0,+\infty)^{2}$, satisfying
\begin{equation*}
P\fp{\rho}{P}+T\fp{\rho}{T}=\rho^{2}\fp{e}{P}\cdot
\end{equation*}
\end{assu}
Introduce the 
$1$-form $\omega$ defined by $T\omega\defn\diff e 
+ P\diff (1/\rho)$, 
where we started using the notation 
$\diff f =(\partial f/\partial T)\diff T+(\partial f/\partial P)\diff P$. 
Assumption~\ref{assu:general} implies that~$\diff \omega=0$. 
Hence, the Poincar\'e's Lemma implies that 
there exists a $C^{\infty}$ function $S=S(P,{T})$, 
defined on $(0,+\infty)^{2}$, 
satisfying the second principle of thermodynamics:
\begin{equation}\label{secondprinciple}
{T} \diff S =  \diff e + P\diff (1/\rho).
\end{equation}
By combining the evolution equations for~$\rho$ and~$e$ 
with \eqref{secondprinciple} written in the form
~$\rho {T}\diff S = \rho \diff e - (p/\rho)\diff \rho$, we get an 
evolution equation for $S$, so that
$$
(\partial_{t}+v\cdot\nabla)\begin{pmatrix}\rho \\ S \end{pmatrix} =
\begin{pmatrix} -\rho & 0 \\ 0 & (\rho {T})^{-1} \end{pmatrix}
\begin{pmatrix} \cn v \\ \kappa\cn(k\nabla T) + Q\end{pmatrix}.
$$
On the other hand, one has
$$
(\partial_{t}+v\cdot\nabla)\begin{pmatrix}\rho \\ S \end{pmatrix}
=J (\partial_{t}+v\cdot\nabla)\begin{pmatrix} P \\ {T} \end{pmatrix}
\quad\text{with}\quad
J\defe
\begin{pmatrix} \partial {\rho}/\partial {P} &
\partial {\rho}/\partial {{T}} \\ 
\partial {S}/\partial {P} & \partial {S}/\partial {{T}} \end{pmatrix}.
$$
Equating both right hand sides and inverting the matrix $J$, we obtain
\begin{equation}\label{computation:1}
\left\{
\begin{aligned}
&(\partial_{t}{P}+v\cdot\nabla P) + a \cn v - \kappa b\cn(k\nabla T)=b Q,\\
&(\partial_{t}T+v\cdot\nabla T) + c\cn v - \kappa d \cn (k\nabla T)=d Q,
\end{aligned}
\right.
\end{equation}
where 
$$
a={\displaystyle \frac{\rho\fpll{S}{{T}}}{\det(J)}}\virgp ~
b={\displaystyle -\frac{\fplll{\rho}{{T}}}{\rho {T}\det(J)}}\virgp ~ 
c ={\displaystyle -\frac{\rho\fpll{S}{P}}{\det(J)}}\virgp ~
d ={\displaystyle \frac{\fplll{\rho}{P}}{\rho {T}\det(J)}}\cdot
$$

To express the coefficients $g_{i}$ and $\chi_{i}$ in terms of
physically relevant quantities, we need some more notations. 
We introduce
\begin{equation}\label{nota:thermo}
\begin{aligned}
K_{{T}} &\defn \frac{1}{\rho}\fp{\rho}{P}\virgp&\qquad
 K_{P} &\defn -\frac{1}{\rho}\fp{\rho}{T}\virgp \qquad
\Rr\defn -\rho\frac{\fplll{S}{P}}{\fplll{\rho}{P}}\virgp\\
C_{P} &\defn T\fp{S}{T}\virgp &\qquad
C_{V} &\defn 
T\frac{\fpll{S}{T}\fpll{\rho}{P}-\fpll{S}{P}\fpll{S}{T}}{\fplll{\rho}{P}}\cdot
\end{aligned}
\end{equation}
The functions~$K_{{T}}$,~$ K_{P}$,~$C_{V}$ and~$C_{P}$ are known 
as the coefficient of isothermal compressibility, 
the coefficient of thermal expansion and the specific heats at 
constant volume and pressure, 
respectively (see Section~$2$ in~\cite{Evans}). 
The function~$\Rr$ generalizes the usual gas constant: 
for perfect gases one can check that~$\Rr={R}$. 

We now have to convert System~\ref{computation:1} into equations for 
the fluctuations $ p$ and $\theta$ as defined by~\eqref{defi:ptheta}. 
Performing a little algebra we find that
\begin{equation*}
\left\{
\begin{aligned}
&\frac{K_{{T}} C_{V}P}{C_{P}}(\partial_{t}p+v\cdot\nabla p)
+\frac{1}{\eps}\cn v 
-\frac{\kappa}{\eps}\frac{ K_{P}}{\rho C_{P}}\cn(kT\nabla \theta)=\frac{1}{\eps}\frac{K_{P}}{\rho C_{P}} Q,\\
&\rho(\partial_{t}v+v\cdot\nabla v)+\frac{1}{\eps}P\nabla p
=\mu\bigl(2\cn(\zeta Dv)+\nabla(\eta\cn v)\bigr),\\
&\rho C_{V}{T}(\partial_{t}\theta+v\cdot\nabla\theta)+
\Rr \rho {T} \cn v -\kappa\cn(kT\nabla \theta)= Q.
\end{aligned}
\right.
\end{equation*}
Hence, $( p, v,\theta)$ satisfies~\eqref{system:NSint} with
\begin{equation}\label{explicitvalues}
g_{1}^{*}\defe \frac{K_{{T}}C_{V}P}{C_{P}}\virgp~ 
g_{2}^{*}\defe \frac{\rho}{P}\virgp~  g_{3}^{*}\defe \frac{C_{V}}{\mathcal{R}}\virgp~
\chi_{1}^{*}\defe \frac{ K_{P}}{\rho C_{P}}\virgp~
\chi_{2}^{*}\defe\frac{1}{P}\virgp~\chi_{3}^{*}\defe \frac{1}{\mathcal{R}\rho{T}}\virgp
\end{equation}
where we used the following notation: for all 
$f\colon (0,+\infty)^{2}\rightarrow \xR$,
\begin{equation}\label{pushforward}
f^{*}(\vartheta,\wp)\defn f(\underline{{T}}e^{\vartheta},\underline{P}e^{\wp}).
\end{equation}

\subsection{Properties of the coefficients}\label{SA2PC}

\begin{assu}\label{assu:general2}
The functions $\rho$ and $e$ are $C^{\infty}$ functions 
of~$(P,{T})\in (0,+\infty)^{2}$ such that, $\rho>0$ and
\begin{equation}\label{AARE} 
\fp{\rho}{P} > 0, \quad \fp{\rho}{{T}} < 0
\quad\mbox{and}\quad
\fp{e}{{T}}\fp{\rho}{P}
>\fp{e}{P}\fp{\rho}{{T}}
\cdot 
\end{equation}
\end{assu}
\begin{remark}
This assumption is satisfied by general equations of state. 
Indeed, \eqref{AARE} just means that the coefficients $K_{T}$, 
$ K_{P}$ and $C_{V}$ are positive.
\end{remark}

The following result prove that Assumptions $\ref{assu:general}$ 
and $\ref{assu:general2}$ imply that our main structural assumption is satisfied.
\begin{proposition}
If Assumptions $\ref{assu:general}$ and $\ref{assu:general2}$ are satisfied,
then $\chi_{1}<\chi_{3}$ and $g_{i}$, $\chi_{i}$ ($i=1,2,3$) are $C^{\infty}$ positive functions.
\end{proposition}
\begin{proof}[Proof]
In view of~\eqref{explicitvalues}, the proof reduces to establishing that
$$
0<K_{{T}},\quad 0< K_{P}, \quad
0< C_{V} < C_{P} \quad\mbox{and}\quad 0 < \Rr   < \frac{C_{P}}{{T} K_{P}}\cdot
$$
The first two inequalities follow from the definitions 
of~$K_{{T}}$ and~$ K_{P}$. To prove the last two, we first establish the Maxwell's identity 
~$\fplll{S}{P}=\rho^{-2}\fpll{\rho}{{T}}$. To see this, by~\eqref{secondprinciple}, 
we compute
\begin{align*}
\fp{S}{P} \diff {T}\wedge \diff P = \diff\bigl( {T} \diff S\bigr) 
= \diff\Bigl\{\diff e + P\diff\Bigl(\frac{1}{\rho}\Bigr)\Bigr\}
= -\frac{1}{\rho^{2}}\fp{\rho}{{T}} \diff P\wedge \diff {T}.
\end{align*}

Since $\fplll{\rho}{{T}}<0$, the Maxwell's identity implies that $\fplll{S}{P}<0$. 
By combining this inequality with $\fplll{\rho}{P}>0$, we find $\Rr>0$. 
Also, the identity~$\fplll{S}{P}=\rho^{-2}\fpll{\rho}{{T}}$ implies that
\begin{equation*}
\frac{C_{P}}{C_{V}}
=\frac{\fpll{S}{{T}}\fpll{\rho}{P}}{\fpll{S}{{T}}\fpll{\rho}{P}-
\rho^{-2}\bigl(\fplll{\rho}{{T}}\bigr)^{2}}\virgp 
\end{equation*}
which proves~$C_{V}<C_{P}$. 

In view of~\eqref{secondprinciple}, the  
assumption 
${\displaystyle \fp{e}{{T}}\fp{\rho}{P}>\fp{e}{P}\fp{\rho}{{T}}}$ 
is equivalent to 
$$
{\displaystyle \fp{S}{{T}}\fp{\rho}{P}>\fp{S}{P}\fp{\rho}{{T}}}\cdot
$$
This inequality has two consequences. Firstly, it implies that
$C_{V}>0$. Secondly, it yields
$$
\frac{{T} K_{P}\Rr}{C_{P}} = \frac{\fpll{S}{P}\fpll{\rho}{{T}}}{\fpll{S}{{T}}\fpll{\rho}{P}}<1.
$$
This concludes the proof.
\end{proof}

We now discuss the physical meaning of the functions~$\Di$ and $\DII$
introduced in $\S$\ref{assu:assus}. 
These are compatibility
conditions between the singular terms and the viscous terms. To see
this, suppose $(p,v,\theta)$ is a smooth solution of~\eqref{system:NSint} and 
let $\Psi=\Psi(\vartheta,\wp)\in C^{\infty}(\xR^{2})$. Then $\psi\defn\Psi(\theta,\eps p)$ satisfies
\begin{multline*}
g_{1}g_{3}\bigl(\partial_{t}\psi+ v\cdot\nabla\psi\bigr)\\
+\Bigl(\underset{=:\Gamma_{1}(\Psi)}{\underbrace{g_{1}\fp{\Psi}{\vartheta}+g_{3}\fp{\Psi}{\wp}}}\Bigr)\cn v
=\kappa\Bigl(\underset{=:\Gamma_{2}(\Psi)}{\underbrace{g_{1}\chi_{3}\fp{\Psi}{\vartheta}+g_{3}\chi_{1}\fp{\Psi}{\wp}}}\Bigr)
\bigl(\cn( k(\theta)\nabla\theta)+Q\bigr),
\end{multline*}
where the coefficients $g_{i}$, $\chi_{i}$, $\fplll{\Psi}{\vartheta}$
and $\fplll{\Psi}{\wp}$ are evaluated at $(\theta,\eps p)$. 
We next show that for appropriate function $\Psi$ one can impose
\begin{equation}\label{A1G12}
[\Gamma_{1}(\Psi)=0 \text{  and  }\Gamma_{2}(\Psi)>0] \quad\text{or}\quad 
[\Gamma_{1}(\Psi)>0 \text{  and  }\Gamma_{2}(\Psi)=0].
\end{equation}
\begin{proposition}\label{prop:additionalcomp}
Assume that Assumptions $\ref{assu:general}$ and $\ref{assu:general2}$
are satisfied and use the notation $\eqref{pushforward}$. The functions $S^{*}$ and $\rho^{*}$ satisfy
\begin{equation}\label{defi:Psigchi}
g_{1}\fp{S^{*}}{\vartheta}=-g_{3}\fp{S^{*}}{\wp}>0, \qquad
g_{1}\chi_{3}\fp{\rho^{*}}{\vartheta}=-g_{3}\chi_{1}\fp{\rho^{*}}{\wp}<0.
\end{equation}
\end{proposition}
\begin{remark}The fact that $\Psi=S^{*}$ (or $\Psi=\rho^{*}$) 
satisfies the first (respectively second) set of conditions in \eqref{A1G12} 
now follows from $\chi_{1}<\chi_{3}$.
\end{remark}
\begin{proof}[Proof]
By~\eqref{explicitvalues} and the definitions given in~\eqref{nota:thermo}, one has
\begin{equation}\label{comp:gPsi}
\frac{g_{1}^{*}}{g_{3}^{*}}=-\frac{P\fpll{S}{P}}{{T}\fpll{S}{{T}}}\cdot
\end{equation}
By definition~\eqref{pushforward}, 
$\fplll{f^{*}}{\vartheta}=\bigl[{T}\fpll{f}{{T}}\bigr]^{*}$ and 
$\fplll{f^{*}}{\wp}=\bigl[P\fpll{f}{P}\bigr]^{*}$. 
This proves that $S^{*}$ satisfies the first identity
in~\eqref{defi:Psigchi}. 
Next, we compute
$$
\frac{\chi_{1}^{*}}{\chi_{3}^{*}}=\frac{\fpll{\rho}{{T}}\fpll{S}{P}}{\fpll{\rho}{P}\fpll{S}{{T}}}\cdot
$$
By~\eqref{comp:gPsi}, this yields
$\chi_{1}^{*}g_{3}^{*}P\fpll{\rho}{P}=-\chi_{3}^{*}g_{1}^{*}{T}\fpll{\rho}{{T}}$.
Which proves that $\rho^{*}$ satisfies the second identity in~\eqref{defi:Psigchi}.
\end{proof}
\begin{remark}
Assumption~(H\ref{assus}) in $\S$\ref{assu:assus} requires, in addition, that 
$F=S^{*}$ and $G=\rho^{*}$ define bijections. This means nothing but the 
fact that the thermodynamic state is completely determined by 
$(P,T)$, or $(P,S)$ or $(\rho,T)$. 
\end{remark}

The following result 
contains an example of equation of state such that 
$\chi_{1}$ depends on $\vartheta$.
\begin{proposition}\label{prop:structurepg}
Assume that the gas obeys Mariotte's law: $P={R}\rho{T}$, for some positive constant $R$, and 
$e=e({T})$ satisfies $C_{V}\defn \fplll{e}{T}>0$. 
Then, 
Assumptions $\ref{assu:general}$ and $\ref{assu:general2}$ are satisfied. Moreover, 
$$
\chi_{1}^{*}=R/(({C}_{V}(T)+{R})P),
$$
so that $\chi_{1}(\vartheta,\wp)$ is independent of $\vartheta$ if and only if $C_{V}$ is constant. 
\end{proposition}

\medskip
\noindent \textbf{Acknowledgments.}
I warmly thank Guy M\'etivier for helpful discussions. 
Thanks also to Didier Bresch, Christophe Cheverry, 
Rapha\"{e}l Danchin and David Lannes for stimulating comments about this work.

\end{document}